%% file: dipl.tex
\documentclass[a4paper,12pt,twoside,BCOR2mm,footexclude,headinclude,normalheadings]{scrreprt}
\pdfoutput=1
\usepackage[utf8]{inputenc}
\usepackage[english]{babel}
\usepackage[T1]{fontenc}
\usepackage{mathpazo}
\usepackage{amsmath}
\usepackage{amssymb}
\usepackage{hyperref}
\usepackage{graphicx}
\usepackage{xypic}
\usepackage[thmmarks,amsmath,hyperref]{ntheorem}

\typearea[5mm]{12}

\author{Jan Essert}
\date{Autumn 2006}

\titlehead{\begin{center}Darmstadt University of Technology\\Department of Mathematics\\Functional Analysis\end{center}}
\subject{Diploma Thesis}
\title{Homological stability for certain classical groups}
\publishers{\vspace{3cm}\includegraphics[height=5cm]{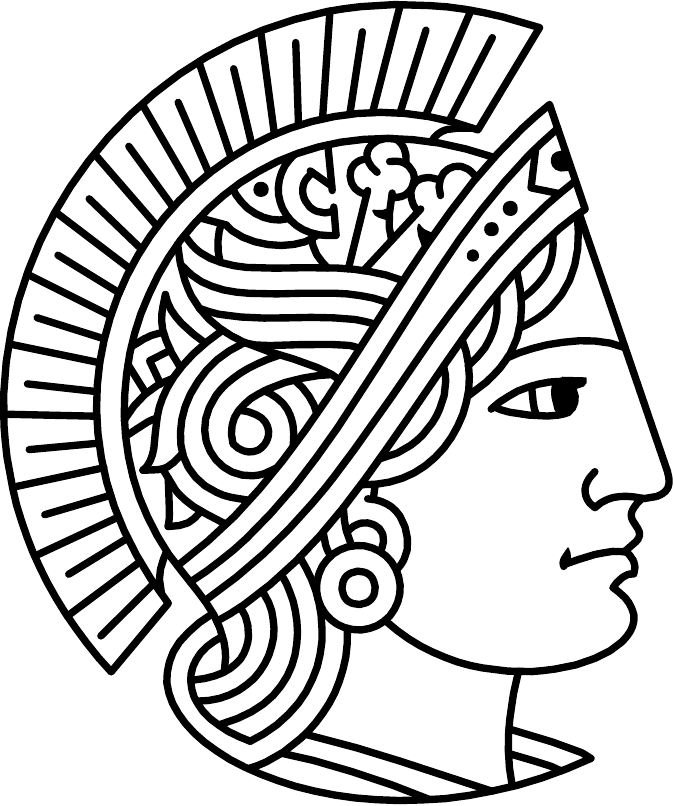}\\[1cm]supervised by Prof.~Dr.~Linus Kramer}

\hypersetup{
pdftitle={Homological stability for certain classical groups},
pdfauthor={Jan Essert},
pdfkeywords={Group Homology, Classical groups, Group theory, General linear groups, Unitary groups, Homological stability}
}

\newcommand{\R}{\mathbb R}
\newcommand{\Q}{\mathbb Q}
\newcommand{\Z}{\mathbb Z}

\newcommand{\K}{\mathbb K}

\newcommand{\C}{\mathbb C}
\newcommand{\F}{\mathbb F}
\newcommand{\A}{\mathbb A}

\newcommand{\bH}{\mathbb H}
\newcommand{\cF}{\mathcal F}
\newcommand{\cC}{\mathcal C}
\newcommand{\cB}{\mathcal B}
\newcommand{\fS}{\mathfrak S}
\newcommand{\Ep}{{'E}}
\newcommand{\Epp}{{''E}}
\newcommand{\cFp}{{'\cF}}
\newcommand{\cFpp}{{''\cF}}
\newcommand{\stern}{*}
\DeclareMathOperator{\id}{id}
\DeclareMathOperator{\im}{im}
\DeclareMathOperator{\spn}{span}
\DeclareMathOperator{\Tor}{Tor}
\DeclareMathOperator{\Ext}{Ext}
\DeclareMathOperator{\Aut}{Aut}
\DeclareMathOperator{\msharp}{\sharp}
\DeclareMathOperator{\Gr}{Gr}
\DeclareMathOperator{\chr}{char}
\DeclareMathOperator{\Gl}{Gl}
\DeclareMathOperator{\Hom}{Hom}

\newcommand{\ruleboxstart}{\vspace{1em}\hrule\nopagebreak}
\newcommand{\ruleboxend}{\nopagebreak\hrule}

\newcommand{\GMG}{\bigl(\begin{smallmatrix} G_1 & M \\ 0 & G_2\end{smallmatrix}\bigr)}
\newcommand{\oMo}{\bigl(\begin{smallmatrix} I_n & M \\ 0 & I_m\end{smallmatrix}\bigr)}
\newcommand{\GoG}{\bigl(\begin{smallmatrix} G_1 & 0 \\ 0 & G_2\end{smallmatrix}\bigr)}

\theoremstyle{plain}
\newtheorem{Definition}{Definition}[chapter]
\newtheorem{Theorem}[Definition]{Theorem}
\newtheorem{Lemma}[Definition]{Lemma}
\newtheorem{Proposition}[Definition]{Proposition}
\newtheorem{Corollary}[Definition]{Corollary}
\newtheorem{Problem}[Definition]{Problem}

\theoremclass{Theorem}
\theoremstyle{break}
\newtheorem{BigTheorem}[Definition]{Theorem}

\theoremstyle{nonumberplain}
\theoremsymbol{\ensuremath\Box}
\theorembodyfont{\rm}
\newtheorem{Proof}{Proof}

\addtokomafont{paragraph}{\rm\bf}
\addtokomafont{subparagraph}{\rm\it}
\begin{document}
\maketitle
\tableofcontents

\vfill
\noindent I would like to thank Sarah and my parents for their help and encouragement during my work on this thesis. In addition, I am very grateful to Linus Kramer for all the advice and support. Many thanks to Ralf Gramlich, Tobias Hartnick, Karsten Hayn, Max Horn, Linus Kramer, Andreas Mars and Nico Weber for proof-reading my drafts and for many valuable corrections, hints and suggestions.

\setcounter{chapter}{-1}

\chapter{Introduction}

This thesis covers homological stability theorems. That means that we consider inclusions of, for example, $\Gl_n\F$ into a point stabilizer of $\Gl_{n+1}\F$. These inclusions induce homomorphisms on homology groups
\[
	H_q(\Gl_n\F)\rightarrow H_q(\Gl_{n+1}\F).
\]
Homological stability is the question: \emph{For which values of $q$ are these induced maps isomorphisms?} We will prove theorems on homological stability for
\begin{itemize}
\item general linear groups over skew-fields with infinite centre and
\item standard unitary groups over $\R$, $\C$ and $\bH$.
\end{itemize}
Both of the proofs are due to Chih-Han Sah, see \cite[section~1 and appendix~B]{Sah:HcL:86}.

\paragraph{} The only prerequisites for this thesis will be some knowledge in homological algebra. The thesis is organised as follows:

In chapter~\ref{c:grouphomology}, we will define group homology, cite standard group homology theorems and prove some introductory results.

Chapter \ref{c:homabgrp} will cover the homology of abelian groups. Here we cite three deep results on the homology of abelian groups and specialise to the homology of vector spaces, which we will need later on in chapter~\ref{c:homisom}.

The main tools for this thesis are spectral sequences. The relevant spectral sequences will be introduced in chapter~\ref{c:specseq}, their existence will not be proven. Most of the chapter is concerned with two spectral sequences associated to a double complex. We investigate these two spectral sequences to get important technical tools for chapter~\ref{c:homstabthms}.

In chapter~\ref{c:homisom}, we present a homological isomorphism theorem by A.~Suslin and prove it.

The heart of this thesis is chapter~\ref{c:homstabthms}, where we state the two resulting theorems. The whole of chapter~\ref{c:homstabthms} is dedicated to the proofs, which will be given in parallel to emphasise the similarities and differences. We finish the thesis with some easy consequences.

\chapter{Homology of groups} \label{c:grouphomology}
In this first chapter, we will briefly introduce group homology. It will be constructed in a very explicit fashion. For a detailed introduction, a more general approach and all the proofs, see \cite{Bro:CoG:82}.

We assume familiarity with the concept of chain complexes, graded modules and homology groups. Knowledge of singular homology may be helpful but is not required. Throughout this chapter, let $G$ be any abstract group. The neutral element in $G$ will always be denoted by $1_G$.

\section{Group rings and modules}

\begin{Definition}
	For any group $G$ we define the \emph{group ring $\Z G$} to be the free abelian group over the elements of $G$ additively. The multiplication is defined by linear extension of the group multiplication.
\end{Definition}
\begin{Definition}
	A left module over $\Z G$ is called a \emph{$G$-module}.
\end{Definition}
	A $G$-module structure can also be seen as an abelian group endowed with a linear left $G$-action. The $G$-module is called \emph{trivial} if this $G$-action is trivial.

\paragraph{Tensor products} Tensor products $M\otimes N$ of modules are only defined if $M$ is a right module and $N$ is a left module over a common ring. We want to form tensor products of $G$-modules over $\Z G$. Note that we can canonically make any left $\Z G$-module $M$ into a right $\Z G$-module by setting
\[
	mg := (g^{-1})m \quad\text{for any}\quad m\in M, g\in G.
\]
Using this construction, we can define tensor products of left $\Z G$-modules $M$ and $N$, denoted by $M\otimes_G N$. In these tensor products we have the following relation:
\[
	(gm\otimes_G gn) = (m(g^{-1})\otimes_G gn) = (m\otimes_G n).
\]
Note that there is some ambiguity here when $M$ already admits a right $\Z G$-module structure, but this will never be a problem in this thesis.

\section{The ordered simplicial chain complex}

Throughout this thesis we will consider chain complexes formed in a particular fashion using a given set $X$ with a left $G$-action.

\begin{Definition}
	For a non-empty set $X$ with a left $G$-action, we define the \emph{ordered simplicial chain complex over $X$}, denoted by $(C_\stern(X),d)$ as follows:
	\[
	C_k(X) := \left\{\begin{array}{ll} 0 & k<0 \\ \Z (X^{k+1}) & k \geq 0 \end{array}\right.,
	\]
	where $\Z(X^{k+1})$ means the free abelian group over $(k+1)$--tuples of elements in $X$ (which will be called \emph{simplices}), endowed with the diagonal $G$-action
	\[
		g(x_0,\ldots,x_k) := (gx_0,\ldots,gx_k),
	\]
	these are $G$-modules. The differentials are defined by
	\[
		d(x_0,\ldots,x_k) := \sum_{j=0}^k(-1)^j(x_0,\ldots, \widehat{x_j},\ldots,x_k)
	\]
	for $k>0$, all other differentials are zero. As usual, $\widehat{x_j}$ indicates omitting this entry.
\end{Definition}
\paragraph{Remark} Note that, although the modules $\Z(X^{k+1})$ are free abelian groups, in general, they are not free as $G$-modules. Representatives of $G$-orbits of simplices are always a generating set, but not necessarily linearly independent. This is only true if the $G$-action is free on $X$.

\begin{Definition}
	For any $x\in X$ and any $k\geq 0$ we define a map, called the \emph{join with $x$} by
\[
	(x\msharp -): C_k(X) \rightarrow C_{k+1}(X),\quad x\msharp(x_0,\ldots,x_k):=(x,x_0,\ldots,x_k)
\]
	and linear extension. This is a homomorphism of abelian groups, and not, in general, of $G$-modules.\label{d:msharp}
\end{Definition}
\begin{Lemma}
	The ordered simplicial chain complex $(C_\stern(X),d)$ is always acyclic. The map $\varepsilon_C$ defined by
	\[
		\varepsilon_C: C_0(X) \rightarrow \Z, \quad (x)\mapsto 1
	\]
	and linear extension induces an isomorphism $H_0(C)\rightarrow \Z$.\label{l:acycl}
\end{Lemma}
\begin{Proof}
	Let $c=(x_0,\ldots,x_k)\in C_k$ for $k> 0$ be a simplex. Pick an arbitrary $x\in X$. Then we know that
\begin{eqnarray*}
	d(x\msharp c) = d(x,x_0,\ldots,x_k) &=& (x_0,\dots,x_k) - \sum_{j=0}^k(-1)^j(x,x_0,\ldots,\widehat{x_j}\,\ldots,x_k) \\
	&=& c - x \msharp (dc)
\end{eqnarray*}
	So, if $dc=0$, i.e.~$c$ is a cycle, then $d(x\msharp c)=c$, which means that $c$ is a boundary.	For $k=0$, we know that
	\[
		H_0(C)=\Z X/\langle x-y | x,y\in X \rangle\cong \Z,
	\]
	and $\ker \varepsilon_C = \langle x-y | x,y\in X \rangle$. Hence $\varepsilon_C$ induces the desired isomorphism.
\end{Proof}

\section{Free resolutions}

\begin{Definition}
	A \emph{free resolution of $\Z$ over $\Z G$} is an exact sequence of the form
\[
	\cdots \rightarrow F_1 \stackrel{d}{\rightarrow} F_0 \stackrel{\varepsilon}{\rightarrow} \Z \rightarrow 0,
\]
	where all $F_i$, $i\geq 0$ are free $G$-modules and where $\Z$ is considered as a trivial $G$-module.
\end{Definition}

\paragraph{Construction} Clearly, $G$ acts on itself via left multiplication. So we form the ordered simplicial chain complex $F_k(G) = C_k(G)$. In this case the $G$-modules $C_k(G)$ are free --- a basis may be formed by those elements having $1_G$ as first entry, which is easily verified.

Since we know that the ordered chain complex is acyclic, the following sequence is exact.
\[
\cdots \rightarrow F_1(G) \stackrel{d}{\rightarrow} F_0(G) \stackrel{\varepsilon_{F(G)}}{\rightarrow} \Z \rightarrow 0,\qquad
\]
\begin{Definition}
	This free resolution $F(G)$ is called the \emph{standard resolution of $\Z$ over $\Z G$}.
\end{Definition}

\section{The definition of group homology}

Now for any $G$-module $M$, we define a new chain complex $F(G,M)$ via
\[
	F_k(G,M):=F_k(G)\otimes_G M,\qquad d_{F(G,M)} = d_F\otimes_G \id.
\]
Since forming a tensor product with $M$ is, in general, not exact, the tensored complex admits non-trivial homology groups.

For $M \in \{\F_p,\Z,\Q\}$, we will always consider $M$ as a trivial $G$-module. Note that for $M=\Z$, a tensor product over $\Z G$ is not a trivial operation. Forming a tensor product $-\otimes_G M$ is equivalent to `dividing out the $G$-action', since for example in $F \otimes_G \Z$ we have
\[
	gf \otimes_G 1_G  = f(g^{-1}) \otimes_G 1_G = f \otimes_G (g^{-1}\cdot 1_G) = f \otimes_G 1_G.
\]
\begin{Definition}
	The homology groups of this complex
	\[
		H_\stern(G,M):=H_\stern(F_\stern(G,M))
	\]
	are called the \emph{homology groups of $G$ with coefficients in $M$}.
\end{Definition}

\paragraph{Cohomology} By using $\Hom_G(-,M)$ instead of $(-\otimes_G M)$, group cohomology can be defined analogously.

\paragraph{Topological interpretation} Group homology with coefficients in a trivial $G$-mod\-ule can also be interpreted as the standard singular homology of a certain classifying topological space associated to the group $G$, the \emph{Eilenberg-MacLane space of type $K(G,1)$}.

In general, the Eilenberg-MacLane spaces of type $K(G,n)$ are certain CW-com\-plex\-es whose $n$-th homotopy group is isomorphic to $G$, while all other homotopy groups vanish. For a construction, see for example \cite[first talk]{Car:AEM:56}.

Using these spaces, a sort of `higher homology groups', called \emph{Eilenberg-MacLane groups} $H_\stern(G,n;M)$ can also be defined. These are the singular homology groups of the Eilenberg-MacLane spaces $K(G,n)$ with coefficients in $M$. 

As stated above, for $n=1$ we are in the situation of {`}standard' group homology, $H_\stern(G,1;M)= H_\stern(G,M)$. The Eilenberg-MacLane groups will not be needed in this thesis. We will, however, invoke a result for these groups in chapter~\ref{c:homabgrp} and specialise to the case $n=1$.

\paragraph{Functoriality} For any $k\geq 0$, we see that $H_k(G)$ is a covariant functor from the category of groups to the category of abelian groups. Homomorphisms $G\rightarrow H$ induce chain maps $F_\stern(G)\rightarrow F_\stern(H)$ which themselves induce homomorphisms of the homology groups.

\paragraph{Remark} At first glance, the above construction may seem somewhat arbitrary. But it can be shown that $F_\stern(G)$ can be replaced by any free, in fact even any projective, resolution of $\Z$ over $\Z G$ without changing $H_\stern(G,M)$. This makes the construction of homology groups appear more natural.

The construction is also justified by the fact that homology (and especially cohomology) groups of small index have quite important algebraic interpretations.

\paragraph{Homology groups of small index} A simple calculation shows that 
\begin{equation}
H_0(G,M)\cong \Z \otimes_G M \cong M/\langle gm - m | g\in G\rangle,\label{e:calch0}
\end{equation}
which especially means $H_0(G)\cong \Z$. It is not much more difficult to show that 
\[
	H_1(G)\cong (G)_{ab} = G/[G,G].
\]

\paragraph{Hopf's theorem} For the second homology group $H_2$, there is an interesting application in topology. For any connected CW-complex $X$ with fundamental group $\pi_1(X)$, there is an exact sequence
\[
	\pi_2(X) \rightarrow H_2(X) \rightarrow H_2(\pi_1(X)) \rightarrow 0,
\]
see for example \cite[Theorem~II.5.2]{Bro:CoG:82}.

\paragraph{$\mathbf H^2$ and group extensions} For any group $G$, the cohomology group $H^2$ is known to classify group extensions with abelian kernel in the following sense:

Let $A$ be an abelian group that is a $G$-module. Take any \emph{group extension of $G$ by $A$}, that is an exact sequence
\[
	0 \rightarrow A \rightarrow E \rightarrow G \rightarrow 0.
\]
The conjugation action of $E$ on its subgroup $A$ induces an action of $E/A\cong G$ on $A$, since $A$ is abelian. We require that this action coincides with the $G$-module structure on $A$.

Two such extensions with groups $E$ and $E'$ are said to be \emph{equivalent} if there is an isomorphism $E\rightarrow E'$ such that the following diagram commutes:
\[
\xymatrix{
&&E\ar[dd]\ar[dr]&&\\
0\ar[r]&A\ar[ur]\ar[dr]&&G\ar[r]&0 \\
&&E'\ar[ur]&&
}
\]
There is a bijection between the equivalence classes of such extensions of $G$ by $A$ and the cohomology module $H^2(G,A)$, see \cite[Theorem~IV.3.12]{Bro:CoG:82}.

\section{Results from group homology}

At this point, we have to cite several results from group homology. Most of them are fundamental for all the arguments later on. We will not give most of the proofs, the missing proofs can be found in \cite{Bro:CoG:82}.\nopagebreak
\paragraph{Additivity} Since the tensor product is an additive functor, we know that
\[
	F_k(G)\otimes_G\bigl(\bigoplus_{j\in J} M_j\bigr) \cong \bigoplus_{j\in J} \bigl(F_k(G)\otimes_G M_j\bigr),
\]
hence
\[
	H_k\bigl(G,\bigoplus_{j\in J} M_j\bigr) \cong \bigoplus_{j\in J}H_k(G,M_j).
\]

\begin{Lemma}[Shapiro]
\label{l:shapiro} Let $G$ be a group, $H \leq G$ a subgroup and $M$ an $H$-module. Then for all $k\geq 0$ we have
\begin{align*}
H_k(H,M) & \cong  H_k(G, \Z G \otimes_H M ),\\
\intertext{the isomorphism being induced by}
( h \otimes_H m )  & \mapsto   (h \otimes_G (1_G \otimes_H m)).
\end{align*}
For $M=\Z$ we have
$$
H_k(H) \cong H_k(G,\Z G \otimes_H \Z)\cong H_k(G,\Z[G/H]).
$$
\end{Lemma}
\begin{Proof}
	See \cite[Proposition~III.6.2]{Bro:CoG:82}.
\end{Proof}
This lemma will mostly be used when considering $M=\Z$. In this case it provides a method to link homology groups with coefficients in any module to homology groups of stabilizers, as we will see later on.

\begin{Proposition}\label{p:conj}
Let $G$ be a group and $g\in G$ be fixed. On $F_k(G)$, let $c_g$ denote the diagonal conjugation by $g$. Then the map
\begin{eqnarray*}
F_k(G,M) &\rightarrow& F_k(G,M)\\
f\otimes_G m &\mapsto& c_g(f)\otimes_G gm
\end{eqnarray*}
induces the identity map on $H_k(G,M)$ for any $k\geq 0$.

\noindent For an abelian group $G$ this specialises to the map
\[
	f \otimes_G m \mapsto f \otimes_G gm
\]
inducing the identity map on homology.
\end{Proposition}
\begin{Corollary}
	Let $G$ be a group, $H\trianglelefteq G$ a normal subgroup and $M$ a $G$-module. Then there is a $G/H$-action on $H_k(H,M)$ induced by
	\begin{eqnarray*}
		G/H \times F_k(H,M) &\rightarrow& F_k(H,M)\\
		(gH,f\otimes_H m) &\mapsto& c_g(f)\otimes_H gm
	\end{eqnarray*}
	for any $k\geq 0$.\label{c:ghaction}
\end{Corollary}
\begin{Proof}
	See \cite[Proposition~III.8.1]{Bro:CoG:82} for a proof of the proposition. Since now any element of $H$ induces the identity on homology, the action above is well-defined.
\end{Proof}
The proposition yields the following result by A.~Suslin (see \cite[Lemma~1.5]{Sus:HoG:84}) as an interesting consequence.

\begin{Lemma} Suppose that $A$ is a commutative ring, $G$ is an abelian group and 
\[
	\varphi: G\rightarrow A^\star
\]
is a group homomorphism (i.e.~$A$ is a $G$-algebra via the homomorphism $\varphi$). 
Then on any $A$-module $M$ the homomorphism $\varphi$ induces a natural $G$-module structure. If $H_0(G,A)=0$, then for any $k\geq 0$ we have $H_k(G,M)=0$.\label{l:vanishinghomology}
\end{Lemma}
\begin{Proof}
	We know by equation \eqref{e:calch0} that $H_0(G,A)=A/I$, where $I$ is the submodule generated by all elements of the form $ga-a$ for all $a\in A$. Since $A$ is a commutative $G$-algebra, the submodule $I$ can be written as the ideal generated as follows
	\[
		I=\langle g\cdot1_A - 1_A | g\in G\rangle=\langle \varphi(g) - 1_A | g\in G\rangle.
	\]
	Since $M$ is an $A$-module and $G$ acts on $M$ by $A$-linear transformations, the groups $H_k(G,M)$ are naturally $A$-modules, the scalar multiplication is induced by
	\[
		a\cdot(f\otimes_G m):=(f\otimes_G a\cdot m).
	\]
	The group $G$ is abelian, so the map
	\begin{equation}
		(f\otimes_G m) \mapsto (f \otimes_G gm) = (f\otimes_G \varphi(g)m)\label{e:multizero}
	\end{equation}
	induces the identity map on the homology modules $H_k(G,M)$ by Proposition~\ref{p:conj}. Hence the map
	\begin{equation*}
		(f\otimes_G m) \mapsto (f\otimes_G (\varphi(g) - 1_A)m)
	\end{equation*}
	induces the zero homomorphism on homology. This means that each $H_k(G,M)$ is annihilated by $I$.

	Now, since $H_0(G,A)=0$, we know that $A=I$. Then each $H_k(G,M)$ is annihilated by $A$, hence by $\varphi(G)$. From this it follows that $H_k(G,M)=0$, since $G$-multiplication as in equation \eqref{e:multizero} must induce the identity map by Proposition~\ref{p:conj}.
\end{Proof}

\pagebreak
\section{Homology with coefficients in a chain complex}

For the construction of spectral sequences and for the proof method used in later chapters, we have to introduce a somewhat technical concept.\nopagebreak
\begin{Definition}
	The \emph{tensor product} of two chain complexes of $G$-modules $(F_\stern,d_F)$ and $(C_\stern,d_C)$ is defined as follows:
	\begin{eqnarray*}
		(F \otimes C)_k &=& \bigoplus_{p+q=k} F_p \otimes_G C_q \\
		d_{F\otimes C}(f\otimes c) &=& d_F f \otimes_G c + (-1)^{\deg f} f \otimes_G d_C c
	\end{eqnarray*}
\end{Definition}

\begin{Definition}
	For any non-negative chain complex of $G$-modules $(C_\stern,d_C)$ we define the \emph{group homology of $G$ with coefficients in $C$} to be the homology of the tensor product complex $(F_\stern(G) \otimes C_\stern)$, that is
	\[
		H_\stern(G,C) = H_\stern(F_\stern(G) \otimes C_\stern).
	\]
\end{Definition}

\noindent This is, of course, a generalisation of group homology with coefficients in a $G$-mo\-dule $M$, where $C_0=M$ and $C_i=0$ otherwise.

We will not need these homology groups themselves, but they will serve as a computational tool to calculate `standard' group homology, as we will see in chapter~\ref{c:specseq}.

\chapter{The homology of abelian groups} \label{c:homabgrp}
In this chapter, we will cite three general theorems on group homology and Ei\-len\-berg-Mac\-Lane groups of abelian groups and look at special cases which will be required in chapter~\ref{c:homisom}.

\section{Graded associative algebras}
The proofs, which we will not give here, use a certain product structure on the graded homology module of abelian groups with coefficients in a commutative ring, the \emph{Pontryagin product}. Therefore we will need some definitions for graded algebras. Throughout this section, let $R$ be any commutative ring. All tensor products of $R$-modules are, of course, tensor products over $R$.

\begin{Definition}
A \emph{graded associative $R$-algebra} is a graded $R$-module $A=\bigoplus_{i=0}^\infty A^i$ with an associative algebra structure which is compatible with the grading, i.e.~
\[
	A^i \cdot A^j \subseteq A^{i+j}.
\]
An element $a\in A^i$ is said to be \emph{homogeneous of degree $i$}.
\end{Definition}
\begin{Definition}
	For two graded $R$-algebras $A$ and $B$ their \emph{tensor product $A\otimes B$} is again a graded $R$-algebra. On $A\otimes B$ we define a multiplication as follows: For homogeneous elements $a,a'\in A$ and $b,b'\in B$, we set
	\[
		 (a\otimes b) \cdot (a'\otimes b') = (-1)^{\deg a' \deg b}(aa'\otimes bb')
	\]
	and extend linearly.
\end{Definition}

\paragraph{Examples} As in singular cohomology, $H^\star(G,R)$ can be endowed with the cup-product and hence can be seen as a graded $R$-algebra.

More important for us, as we will see in the next section, for any \emph{abelian} group $G$ and any commutative ring $R$, the graded homology module $H_\star(G,R)$ can be endowed with the aforementioned Pontryagin product, which makes it into a graded associative algebra.

\subsection{The Pontryagin product \texorpdfstring{---}{-} the algebra \texorpdfstring{$H_\stern(G,R)$}{H(G,R)}}

\begin{Definition}[The homology cross-product]
	Let $G$ and $G'$ be groups and let $F$ and $F'$ be free resolutions of $\Z$ over $\Z G$ and $\Z G'$, respectively. For any $G$-module $M$ and any $G'$-module $M'$, there is a natural map
	\[
		(F\otimes_G M)\otimes(F'\otimes_{G'} M') \rightarrow (F\otimes F')\otimes_{G\times G'} (M\otimes M').
	\]
	It is not difficult to see that this map induces a homology map
	\[
	\times: H_p(G,M)\otimes H_q(G',M') \rightarrow H_{p+q}(G\times G',M\otimes M'),
	\]
	called the \emph{homology cross-product}.
\end{Definition}

\noindent Now let $G$ be an abelian group and let $R$ be a commutative ring, viewed as a trivial $G$-module. For an abelian group, the multiplication $G\times G \rightarrow G$ is a group homomorphism. Hence the map
\[
	\mu: (G\times G, R\otimes R) \rightarrow (G,R),\quad\mu((g,g'),(r\otimes r')) = (gg',rr')
\]
induces a homomorphism $\mu_\stern$ on homology groups.
\begin{Definition}[The Pontryagin product]
	The composite map
	\[
		H_p(G,R)\otimes H_q(G,R) \stackrel{\times}{\rightarrow} H_{p+q}(G\times G,R\otimes R) \stackrel{\mu_\stern}{\rightarrow} H_{p+q}(G,R)
	\]
	is called the \emph{Pontryagin product} on $H_\stern(G,R)$. It can be shown that, with this product, $H_\stern(G,R)$ becomes a graded, associative and anti-commutative algebra, see for example \cite[section~V.5]{Bro:CoG:82}.
\end{Definition}

\subsection{The tensor algebra}
Now let $M$ be any $R$-module.
\begin{Definition}
	The \emph{tensor algebra over $M$} is the direct sum
\[
	T(M) := \bigoplus_{k=0}^\infty T^k(M),
\]
	where $T^k(M):= \overbrace{M \otimes \cdots \otimes M}^{k\text{ times}}$ and $T^0(M):=R$. Since $T^k(M)\otimes T^j(M)$ is canonically isomorphic to $T^{k+j}(M)$, we define the product of two homogeneous elements as their tensor product and extend linearly.
\end{Definition}
\paragraph{Remark}
The tensor algebra is an associative, unital and graded algebra that has the following universal property:

Given any unital, associative $R$-algebra $A$, any $R$-module $M$ and any module homomorphism $M\rightarrow A$, there is a canonical algebra homomorphism $T(M)\rightarrow A$ such that the following diagram commutes:
\[
	\xymatrix{
	M\ar[d]\ar[r] & A \\
	T(M)\ar@.[ru]
	}
\]
\paragraph{Action of the symmetric group} The symmetric group $S_k$ acts on $T^k(M)$ as follows: For any $\sigma\in S_k$ we set
\[
	(v_1\otimes\cdots\otimes v_k)^\sigma := v_{\sigma(1)}\otimes\cdots\otimes v_{\sigma(k)}.
\]
\subsection{The exterior algebra}
Starting with the tensor algebra, we construct various algebras by dividing out appropriate ideals.
\begin{Definition}
	The \emph{exterior algebra $\Lambda(M)$} is defined as follows. Set 
	\[
	\Lambda^k(M):=T^k(M)/\mathfrak a^k,
	\]
	where
	\[
		\mathfrak a^k := \bigl\langle v \in T^k(M) : \exists 1\leq i < j \leq k,\, v^{(i,j)} = v \bigr\rangle.
	\]
	Here $(i,j)\in S_k$ is the transposition swapping $i$ and $j$. The multiplication on $T(M)$ induces a multiplication on $\Lambda(M)=\bigoplus \Lambda^k(M)$.
	
	\noindent In fact, $\mathfrak a:=\bigoplus \mathfrak a^k$ is a graded ideal of $T(M)$ and we have $\Lambda(M)=T(M)/\mathfrak a$.
\end{Definition}
\paragraph{Properties} The exterior algebra is an associative, unital, graded and anti-com\-mu\-ta\-tive algebra that has the following universal property:

Given any unital, associative $R$-algebra $A$, any $R$-module $M$ and a module homomorphism $j:M\rightarrow A$ with the property that for any $v\in M$, $j(v)\cdot j(v) = 1$. Then there is a unique algebra homomorphism $f:\Lambda(M)\rightarrow A$ such that the following diagram commutes:
\[
	\xymatrix{
	M\ar[d]\ar[r]^j & A \\
	\Lambda(M)\ar@.[ru]_f
	}
\]
\subsection{The symmetric algebra}
\begin{Definition}
The \emph{symmetric algebra $S(M)$} is constructed as follows: Set
\[
	S^k(M) := T^k(M)/\mathfrak b^k,
\]
where
\[
	\mathfrak b^k := \bigl\langle v^\sigma - v: v\in T^k(M), \sigma\in S^k \bigr\rangle.
\]
The multiplication on $T(M)$ induces a multiplication on $S(M)=\bigoplus S^k(M)$.

\noindent Again, the direct sum $\mathfrak b:=\bigoplus \mathfrak b^k$ is a graded ideal of $T(M)$ and $S(M)=T(M)/\mathfrak b$.
\end{Definition}
\paragraph{Properties} $S(M)$ is an associative, graded, unital and commutative algebra having the following universal property:

Given any unital, associative and commutative $R$-algebra $A$, any $R$-module $M$ and any module homomorphism $M\rightarrow A$, there is a unique algebra homomorphism $S(M)\rightarrow A$ such that the following diagram commutes:
\[
	\xymatrix{
	M\ar[d]\ar[r] & A \\
	S(M)\ar@.[ru]
	}
\]

\paragraph{A basis for $S(M)$} If $M$ is a free $R$-module with basis $\{e_i\}_{i\in I}$, suppose $I$ is a totally ordered set, then $S^k(M)$ has a basis
\[
	\cB(\fS^k(M)) = \bigr\{e_{i_1}\otimes\cdots\otimes e_{i_k} : i_1\leq\ldots\leq i_k \in I \bigr\}.
\]
\subsection{The shuffle algebra}
At last, we construct a perhaps more uncommon algebra. We begin with another product on $T(M)$, the \emph{shuffle product}, defined by linear extension of
\begin{eqnarray*}
	* : T^k(M) \otimes T^j(M) &\rightarrow& T^{k+j} (M) \\
	(v_1\otimes\cdots\otimes v_k, w_{k+1}\otimes\cdots\otimes w_{k+j}) &\mapsto& \sum_{\sigma\in S_{k+j}^{k,j}} (v_1\otimes\cdots\otimes v_k \otimes w_{k+1} \otimes\cdots\otimes w_{k+j})^\sigma,
\end{eqnarray*}
where $S_{k+j}^{k,j}\subseteq S_{k+j}$ are the \emph{shuffles}, i.e.~permutations that leave the order of the $v$ and $w$ invariant, which means
\begin{align*}
	\sigma(i)<\sigma(j)\quad &\text{for }& 1\leq & \,i < j \leq k, \\
	\sigma(i)<\sigma(j)\quad &\text{for }& k+1\leq & \,i < j \leq k+j.
\end{align*}
\begin{Definition}
We consider the submodules $\fS^k(M)\leq T^k(M)$, defined as follows
\[
	\fS^k(M) := \bigl\{ v \in T^k(M) : \forall \sigma \in S_k,\, v^\sigma = v\bigr\}.
\]
The direct sum $\fS(M):=\bigoplus \fS^k(M)$ endowed with the shuffle product is an associative, graded, commutative algebra, called the \emph{shuffle algebra $\fS(M)$}.
\end{Definition}

\paragraph{A basis for $\fS(M)$} If $M$ is a free $R$-module with basis $\{e_i\}_{i\in I}$, again suppose that $I$ is a totally ordered set, then $\fS^k(M)$ admits a natural basis
\[
	\cB(\fS^k(M)) = \bigr\{ e_{i_1\cdots i_k} := \sum_{\sigma\in S_k} (e_{i_1}\otimes\cdots\otimes e_{i_k})^\sigma : i_1\leq\ldots\leq i_k \in I\bigr\}.
\]

\begin{Proposition}
	For any free $R$-module $M$, the symmetric algebra $S(M)$ and the shuffle algebra $\fS(M)$ are isomorphic as graded $R$-modules (\emph{not} as $R$-algebras!).\label{p:symmshuf}
\end{Proposition}
\begin{Proof}
	For all $k$ the map
	\begin{eqnarray*}
		\cB(S^k(M)) &\rightarrow& \cB(\fS^k(M)) \\
		e_{i_1}\otimes\cdots\otimes e_{i_k} &\mapsto& e_{i_1\cdots i_k}
 	\end{eqnarray*}
	is a bijection of the bases and can hence be extended to a module isomorphism $S^k(M)\rightarrow\fS^k(M)$. Hence we can construct a graded module isomorphism $S(M)\rightarrow\fS(M)$.
\end{Proof}

\section{Homology of an abelian group with coefficients in \texorpdfstring{$\F_p$}{Fp} (\texorpdfstring{$p$}{p} an odd prime)}

Let $\F_p$ denote the field of $p$ elements, where $p$ is prime. For the rest of this chapter, abelian groups are written additively.

In this and the next section, we will invoke a result from \cite[talks 8--10]{Car:AEM:56}, where the homology of abelian groups with coefficients in $\F_p$ for any prime $p$ is calculated. The result is much more general than what we need here, in particular, all Ei\-len\-berg-Mac\-Lane groups $H_\star(G,n;\F_p)$ are calculated.

We will state the result here and work our way backwards through the definitions. By specialisation, we arrive at the (much simpler) result for group homology.

Whenever we speak of the homology groups of an abelian group $G$ with coefficients in a field $F$, we consider the action of $G$ on $F$ to be trivial, i.e.~$F$ as a trivial $G$-module.

\begin{BigTheorem}[Cartan]
	For any abelian group $G$, any odd prime $p$ and any $n\geq 1$, we have
	\[
		H_\stern(G,n;\F_p) \cong U(M^{(n)}).
	\]
	This is an isomorphism of $\F_p$-algebras.
\end{BigTheorem}
Of course we need to explain the terminology used above. The graded $\F_p$-vector space $M^{(n)}$ is constructed as a direct sum of $\F_p$-vector spaces $G/pG$ and
\[
	_pG := \{g\in G: pg = 0\}.
\]
The number and respective degrees of these summands is obtained by counting words in a certain word monoid. We will not go into details here, it suffices to say that
\[
	M^{(1)} = \underbrace{G/pG}_{\text{degree }1} \ \ \oplus \underbrace{{_p}G}_{\text{degree }2}.
\]
By construction of $M^{(n)}$, we can pick a basis of homogeneous elements of positive degree. The algebra $U(M^{(n)})$ is defined as the tensor product $\Lambda(M^-)\otimes\fS(M^+)$. Here $M^+$ is the subspace of $M^{(n)}$ generated by the basis elements of even degree, $M^-$ is the subspace generated by the basis elements of odd degree. For $n=1$ the theorem becomes much simpler:

\begin{Corollary}
	For any abelian group $G$ and any odd prime $p$, we have
	\[
		H_\stern(G,\F_p) \cong \Lambda(G/pG) \otimes \fS(_pG),
	\]
	where the generators of $\Lambda(G/pG)$ have degree one and the generators of $\fS({_p}G)$ have degree two.
\end{Corollary}
We will need the following special case.
\begin{Corollary} For any odd prime $p$ and any $\F_p$-vector space $V$ we have
	\[
		H_\stern(V,\F_p) \cong \Lambda(V) \otimes S(V),
	\]
	the generators of $\Lambda(V)$ have degree one and the generators of $S(V)$ have degree two. This is an isomorphism of graded $\F_p$-vector spaces (and not, in general, of graded $\F_p$-algebras).\label{c:homvectp}
\end{Corollary}
\begin{Proof}
	For an $\F_p$-vector space, we have $V/pV = {_p}V = V$, then apply Proposition~\ref{p:symmshuf}.
\end{Proof}

\section{Homology of an abelian group with coefficients in \texorpdfstring{$\F_2$}{F2}}

The corresponding theorem for $\F_2$-coefficients in \cite[talk 10]{Car:AEM:56} reads quite similarly:

\begin{BigTheorem}[Cartan]
	For any abelian group $G$ and any $n\geq 2$, we have
	\[
		H_\stern(G,n;\F_2) \cong \fS(M'^{(n)}).
	\]
	This is an isomorphism of $\F_2$-algebras. If $2G=0$ this is also true for $n=1$.
\end{BigTheorem}

The construction of $M'^{(n)}$ is quite similar to the one we omitted above, it amounts to counting the number of words in a certain word monoid and adding according numbers of $G/2G$ and $_2G$, all of which are $\F_2$-vector spaces.

In this case, since the word monoid is different, we arrive at $M'^{(1)}=G/2G$, graded by the degree $1$. Furthermore, since any abelian group satisfying $2G=0$ is already an $\F_2$-vector space, we have the following, simpler theorem:

\begin{Theorem}
	For any $\F_2$-vector space $V$, there is an isomorphism of graded $\F_2$-algebras
	\[
		H_\stern(V,\F_2) \cong \fS(V).
	\]
\end{Theorem}

\begin{Corollary} \label{c:homvect2}
	For any $\F_2$-vector space $V$, there is an isomorphism of graded $\F_2$-vector spaces
	\[
		H_\stern(V,\F_2) \cong S(V).
	\]
\end{Corollary}

\begin{Proof}
	Apply Proposition~\ref{p:symmshuf}.
\end{Proof}

\newpage
\section{Homology of an abelian group with coefficients in \texorpdfstring{$\Q$}{Q}}

We cite \cite[theorem V.6.4]{Bro:CoG:82}:

\begin{Theorem}
	Assume that $R$ is a principal ideal domain.
	\begin{enumerate}
	\item There is an injective homomorphism of $R$-algebras $\psi:\Lambda(G\otimes R) \rightarrow H_\stern(G,R)$ for every abelian group $G$, which is a split injection if $G$ is finitely generated.
	\item Suppose that every prime $p$ for which $G$ has $p$-torsion is invertible in $R$, then $\psi$ is an isomorphism.
	\item If $R$ has characteristic zero (e.g., if $R=\Z$), then $\psi$ is an isomorphism in dimension $2$.
	\end{enumerate}
\end{Theorem}

In $R=\Q$, of course, every prime is invertible, hence for a $\Q$-vector space we obtain the following corollary.

\begin{Corollary}
	For any $\Q$-vector space $V$, there is an isomorphism of graded $\Q$-algebras
	\[
		H_\stern(V,\Q) \cong \Lambda(V).
	\]\label{c:homvectq}
\end{Corollary}

\chapter{Spectral sequences} \label{c:specseq}
The main computational tools for this thesis will be spectral sequences. While we will not prove the existence of the relevant spectral sequences, we will, however, give a definition of spectral sequences, introduce the required spectral sequences and prove or cite the required properties.

\section{Definition}

We will only define first quadrant spectral sequences. At first glance, the definitions below may seem completely arbitrary, it may not even be clear that spectral sequences exist. We will show afterwards that spectral sequences arise quite naturally. 

\begin{Definition}
	A \emph{first quadrant bicomplex $(C_{p,q},d_{p,q})$ of bidegree $(a,b)$} is a bigraded module with  differentials
	\[
		d_{p,q} : C_{p,q} \rightarrow C_{p+a,q+b}
	\]
	where additionally $C_{p,q}=0$ for $p<0$ or $q<0$.
\end{Definition}

\noindent A bicomplex thus is made of infinitely many chain complexes aligned in a special fashion. We can calculate homology  and arrive at a new bigraded module.

Now, if we had a new differential on the homology bigraded module, we could continue forming homology modules of homology modules. The following definition makes this idea precise.

\begin{Definition}
	A \emph{first quadrant spectral sequence $E$ of homological type} is a sequence of first quadrant bicomplexes $(E^r_{p,q})_{r\geq 0}$ with bidegrees $(-r,(r-1))$, such that
	\[
		E^{r+1} = H(E^r).
	\]
\end{Definition}

\begin{figure}
	\input{firstpages} \caption{The first three pages of a first-quadrant spectral sequence of homological type} \label{f:firstthreepages}
\end{figure}
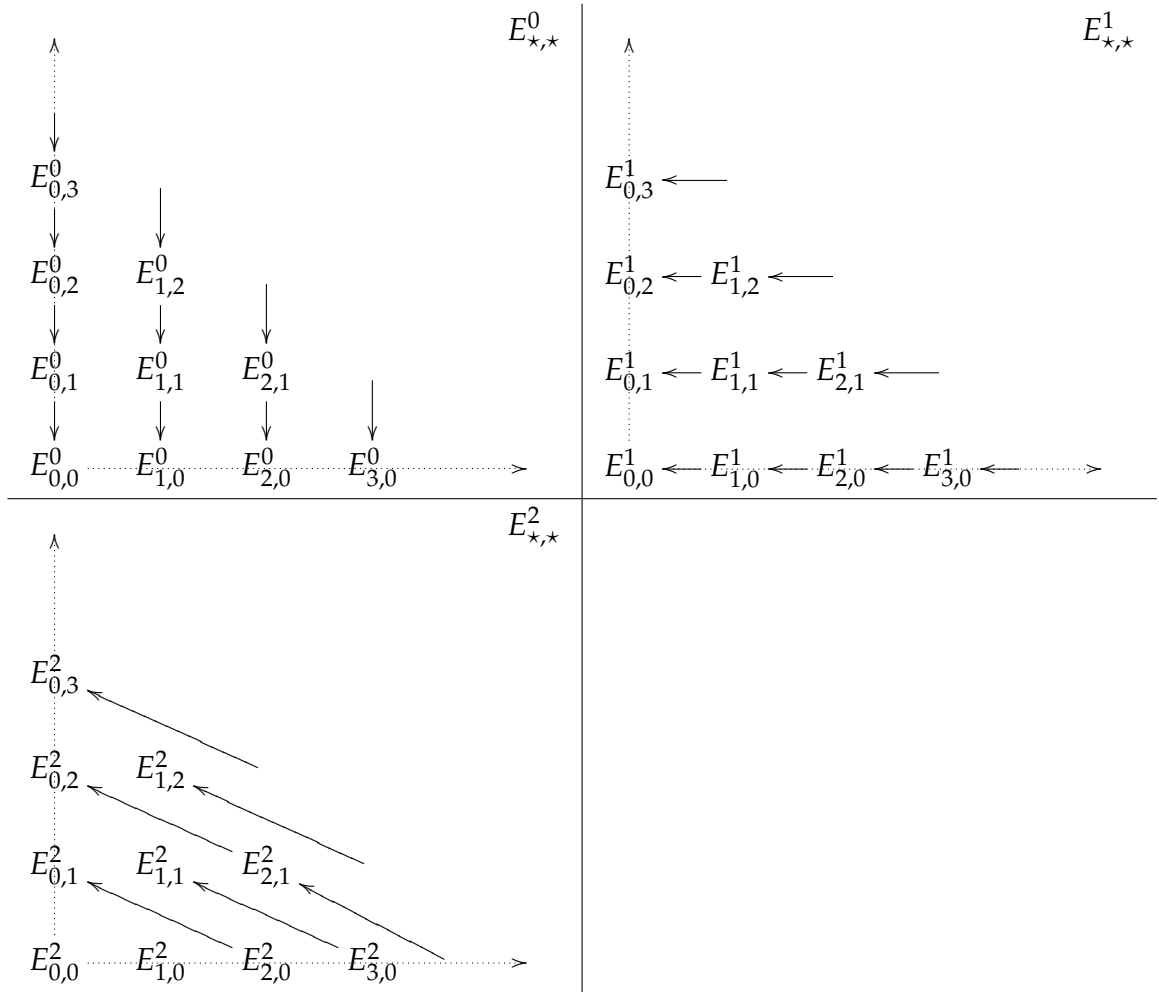

\paragraph{Visualisation} We picture this as a book, whose pages are indexed by $r$. On each page, the graded module $E^r$ can be visualised as the modules $E^r_{p,q}$ sitting on the point $(p,q)$ in the real plane. Turning a page corresponds to calculating homology groups with respect to the differentials on this page. The differential arrows get longer and more slanted as the pages turn. For a sketch of the first three pages of a first-quadrant spectral sequence, see figure \ref{f:firstthreepages}.

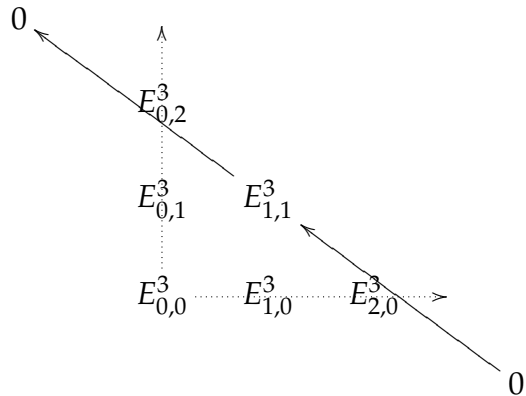
\begin{figure}[htb]
	\input{convergence} \caption{Convergence of a first-quadrant spectral sequence} \label{f:convergence}
\end{figure}

\paragraph{Convergence} Since we only consider first quadrant spectral sequences, for each pair $(p,q)$ all differentials from and into $E^r_{p,q}$ vanish from the $(\max(p,q)+2))$-th page on. The incoming and outgoing arrows begin or end in trivial groups, see figure \ref{f:convergence}. This means that every module $E^r_{p,q}$ stabilizes eventually, we call this stable module $E^\infty_{p,q}$. For any pair $(p,q)$ there is an index $r$, such that
\[
E^r_{p,q}=E^{r+1}_{p,q}=\cdots=:E^\infty_{p,q}.
\]
We call the resulting bigraded module $E^\infty$ the \emph{abutment} of the spectral sequence and we say that the spectral sequence $E$ \emph{converges} to $E^\infty$.
\paragraph{Collapse} If, for any reason, all differentials are zero from a page $E^r$ on, then, obviously
\[
	E^r=E^{r+1}=\cdots=E^\infty,
\]
and we say that the spectral sequence \emph{collapses} on the $r$-th page.

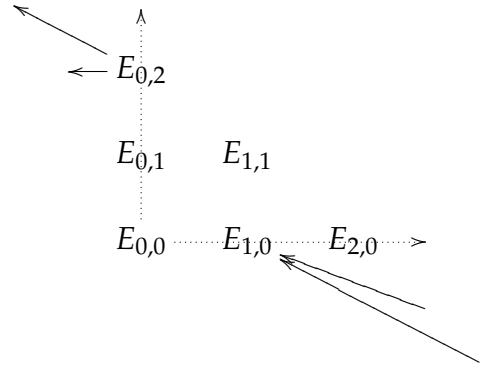
\begin{figure}[hbt]
	\input{edgehom} \caption{Edge homomorphisms: Incoming and outgoing differentials} \label{f:edgehom}
\end{figure}

\paragraph{Edge homomorphisms} For a spectral sequence, let us consider the modules $E^r_{p,0}$ for $r\geq 2$. Since all `incoming' differentials from the second page on are zero (see figure \ref{f:edgehom}), each $E^{r+k}_{p,0}$ is a submodule of $E^r_{p,0}$ for any $k\geq 0$. This of course is also true for $E^\infty_{p,0}$.

So for any $r\geq 2$ and $p\geq 0$ we have a canonical inclusion
\[
	\iota_p: E^\infty_{p,0} \hookrightarrow E^r_{p,0}.
\]

\noindent On the other hand, every `outgoing' differential of the groups $E^r_{0,q}$ is zero from the first page on. Then every $E^{r+k}_{0,q}$ is a quotient of $E^r_{0,q}$ for any $k\geq 0$. So we have a canonical projection
\[
	\pi_q: E^r_{0,q} \twoheadrightarrow E^\infty_{0,q}
\]
for any $r\geq 1$ and any $q\geq 0$. These maps are called \emph{edge homomorphisms}. In the `book picture' they can be imagined as arrows pointing up, respectively down at the `edges of the book'.
\section{The spectral sequence associated to a filtered chain complex}

Spectral sequences arise very naturally in the context of filtered chain complexes. In this section, we take all graded modules (chain complexes) $M_q$ to be non-negative, that means $M_q=0$ for $q<0$.

\begin{Definition}
	A \emph{filtered graded module (chain complex)} is a sequence of graded submodules (subcomplexes) $(\cF^p M)_{p\in\Z}$ of a graded module (chain complex) $(M_q)$ such that $\cF^p M\subseteq \cF^{p+1} M$. We will always use filtrations where $\cF^{-1} M = 0$ and $\cF^p M_p = M_p$, these are called \emph{canonically bounded}.
\end{Definition}

\paragraph{Remark} A filtration on a chain complex $C$ induces a natural filtration on the graded module $H(C)$ as follows:
\[
	\cF^p H(C) := \im\bigl( H(\cF^p C) \rightarrow H(C)\bigr).
\]

\begin{Definition}
	To any filtered graded module $(M_q)$, especially to a filtered chain complex, we associate a bigraded module
\[
	\Gr_{p,q} M := \cF^p M_{p+q} / \cF^{p-1} M_{p+q}.
\]
	Note that for a canonically bounded filtration, $\Gr_{p,q}M=0$ for $p<0$ or $q<0$. So for a chain complex $C$ with a canonically bounded filtration, $\Gr_{p,q}C$ is a first-quadrant bicomplex.
\end{Definition}
\paragraph{Construction}For a filtered chain complex, we have
\[
	E^0_{p,q} := \Gr_{p,q} C = \cF^p C_{p+q} / \cF^{p-1} C_{p+q}
\]
and the differentials on $C$ induce differentials of bidegree $(0,-1)$
\[
	d^0_{p,q}:\Gr_{p,q} C \rightarrow \Gr_{p,q-1}C.
\]

\noindent This means that $\Gr_{p,q}$ splits into `vertical' chain complexes, for which we can calculate homology groups. We arrive at an object called
\[
	E^1_{p,q} := H_q(\Gr_{p,\stern}) = \ker d^0_{p,q} /\im d^0_{p,q+1}.
\]
If we inspect the involved groups more closely, we can see that the differential $d_C$ induces new, `horizontal' differentials (of bidegree $(-1,0)$) on $E^1$. We can again calculate homology groups and arrive at a graded module $E^2$. Here the differential $d_C$ induces new differentials, this time of bidegree $(-2,1)$. We again calculate homology and proceed.

At every step, the differential $d_C$ induces new differentials with the right bidegree. A rather tedious calculation shows that this is a spectral sequence and that it converges to the graded module associated to the homology of $C$,
\[
	E^\infty_{p,q} = \Gr_{p,q} H(C) = \cF^p H_{p+q}(C) / \cF^{p-1} H_{p+q}(C).
\]
This results in the following theorem:

\begin{Theorem}
	For any filtered chain complex $C$, there is a spectral sequence of homological type with first term
	\[
		E^0_{p,q}=\Gr_{p,q} C
	\]
	converging to 
	\[
		E^\infty_{p,q} = \Gr_{p,q} H(C).
	\]
	For a canonically bounded filtration, this is a first-quadrant spectral sequence. The standard notation for convergence is
	\[
		E^0_{p,q}=\Gr_{p,q} C \quad\Rightarrow\quad H_{p+q}(C).
	\] \label{t:exspecseq}
\end{Theorem}
\begin{Proof}
	We will omit this calculation, it can be found in \cite[section~2.2]{McC:UGS:01} or \cite[section~5.4]{Wbl:IHA:94}. A very instructive article to understand where the construction comes from is \cite{Cho:Ych:06}.
\end{Proof}

\paragraph{Convergence} Note that, although the notation suggests that the abutment of this spectral sequence is $H_\stern(C)$, we actually only get something weaker, namely just $\Gr_{\stern,\stern}H(C)$. In general, it may not be possible to reconstruct $H_\stern(C)$ from $\Gr_{\stern,\stern}H(C)$, this is an extension problem. We will only deal with cases where this is possible, however.

\section{The spectral sequences associated to a double complex}

For us, the interesting case of a filtered chain complex will be the following: Suppose we have two non-negative chain complexes over the same commutative ring $R$, called $F$ and $C$. We form the tensor product $(F\otimes C)$ and we want to calculate $H(F\otimes C)$.
\paragraph{Construction} There are two canonical filtrations on $(F\otimes C)$:
	\begin{eqnarray*}
		\cFp^n (F\otimes C)_k &:=& \bigoplus_{\substack{p+q=k\\ p\leq n}} F_p\otimes C_q \\
		\cFpp^n (F\otimes C)_k &:=& \bigoplus_{\substack{p+q=k\\ q\leq n}} F_p\otimes C_q.
	\end{eqnarray*}
Then, of course, by theorem~\ref{t:exspecseq} we know that there are two first-quadrant spectral sequences, accordingly called $\Ep$ and $\Epp$, both converging to $H_\stern(F\otimes C)$. Note that the abutments $\Ep^\infty$ and $\Epp^\infty$ are different, though, since they correspond to different filtrations.

\paragraph{First pages} By examining the graded modules $\Gr_{\stern,\stern} (F\otimes C)$ for these two filtrations, one can easily see that
\[
	\Ep^0_{p,q} = F_p \otimes C_q\qquad\text{and}\qquad \Epp^0_{p,q} = C_p \otimes F_q.
\]
The differentials on the first page are $\id_F\otimes d_C$ and, respectively, $\id_C \otimes d_F$. By forming homology groups for this `vertical' differential, we see that we have different coefficient modules for every column. We arrive at
\[
	\Ep^1_{p,q} = H_q(C,F_q)\qquad\text{and}\qquad \Epp^1_{p,q} = H_q(F,C_p).
\]
By using the explicit construction of these spectral sequences (which we omitted in this thesis), one can see that the differentials on these pages are induced by $d_F\otimes \id_C$ and $d_C \otimes \id_F$. Thus we have the following

\begin{Theorem}
	For any two non-negative chain complexes $F$ and $C$ there are two first-qua\-drant spectral sequences
	\[
	\begin{array}{cc} \Ep^1_{p,q} &= H_q(C,F_p) \\ \Epp^1_{p,q} &= H_q(F,C_p) \end{array} \quad\Rightarrow\quad H_{p+q}(F\otimes C).
	\]
\end{Theorem}

\noindent We inspect the first spectral sequence more closely.

\begin{Proposition}
	For a free non-negative chain complex $F$ and an acyclic non-negative chain complex $C$, we have
	\[
	 \cFp^{p-1}H_p(F \otimes C) = 0.
	\]
	Hence every element in $H_p(F\otimes C)$ has a representative element in $F_p\otimes C_0$. There is an isomorphism
	\begin{align*}
		H_p(F \otimes C) &\rightarrow H_p(F) \\
		\intertext{induced by}
		F_p\otimes C_0 &\rightarrow F_p \otimes \Z \\
		f \otimes c &\mapsto f \otimes \varepsilon_C(c)
	\end{align*} \label{p:homologyofg}
\end{Proposition}

\begin{Proof} Since forming a tensor product with the free modules $F_p$ is exact, we know that
	\[
		\Ep^1_{p,q} = H_q(C,F_p) \cong H_q(C)\otimes F_p = \left\{\begin{array}{cc}H_0(C)\otimes F_p & q=0 \\ 0 & q\neq 0 \end{array}\right..
	\]
We can calculate `horizontal' homology groups and arrive at
	\[
		\Ep^2_{p,q} = \left\{\begin{array}{cc}H_p(F,H_0(C)) & q=0 \\ 0 & q\neq 0 \end{array}\right..
	\]
Then, obviously, $\Ep$ collapses on the second page. This means that
	\[
		\Ep^\infty_{p,q} = \left\{\begin{array}{cc}H_p(F,H_0(C)) & q=0 \\ 0 & q\neq 0 \end{array}\right..
	\]
This is a case where we can reconstruct $H_\stern(F\otimes C)$, we know that
\[
	\Gr_{p,q} H(F\otimes C) = \cFp^p H_{p+q}(F\otimes C) / \cFp^{p-1} H_{p+q}(F\otimes C),
\]
especially
\[
	\Gr_{p,0} H(F\otimes C) = \cFp^p H_{p}(F\otimes C) / \cFp^{p-1} H_{p}(F\otimes C).
\]
Using the above information, we know that the $\cFp$-filtration is trivial in the sense that
\[
	\cFp^{p-1}H_p(F\otimes C) = 0
\]
and hence
\[
	H_p(F\otimes C) = \cFp^p H_p(F\otimes C)/\cFp^{p-1} H_p(F\otimes C) = H_p(F,H_0(C)) \cong H_p(F).
\]
Since $\varepsilon_C$ induces an isomorphism $H_0(C)\cong\Z$ by Lemma~\ref{l:acycl}, we get the desired isomorphism.
\end{Proof}

\paragraph{Application to group homology} Let $G$ be a group. If we take $F=F(G)$, the standard resolution of $G$, and pick $C$ to be an acyclic non-negative complex of $G$-modules, we have two spectral sequences converging to $H_\stern(G,C)$. By the above proposition we know that $H_\stern(G,C)$ is isomorphic to $H(G)$. But we still have the second spectral sequence! We proved the following Corollary which will be a powerful tool to calculate group homology.

\begin{Corollary}
	For any group $G$ and any acyclic non-negative chain complex $C$ of $G$-mo\-dules, there is a first-qua\-drant spectral sequence
	\[
		\Epp^1_{p,q} = H_q(G,C_p) \quad\Rightarrow\quad H_{p+q}(G,C) \cong H_{p+q}(G).
	\]
	The differential on the first page ${''d}^1$ is induced by $\id_{F(G)} \otimes_G d_C$.
	
\noindent Every element in $H_p(G,C)$ has a representative element of the form $f\otimes_G c\in F_p\otimes_G C_0$ and the isomorphism 
\begin{align*}
	H_p(G,C) &\rightarrow H_p(G)\\
	\intertext{is induced by}
		F_p(G)\otimes C_0 &\rightarrow F_p(G)\otimes \Z \\
		(f\otimes_G c) &\mapsto (f\otimes_G \varepsilon_C(c)).
\end{align*} \label{c:acyclicC}
\end{Corollary}
\paragraph{Remark} The previous Corollary is the only reason we introduced homology with coefficients in a chain complex. The interesting part is only the spectral sequence, which gives us much structure to use in the final proofs.
\paragraph{Edge homomorphisms} For this last spectral sequence, there are, of course, the canonical edge homomorphisms
\[
	\pi_q: \Epp^1_{0,q} = H_q(G,C_0) \twoheadrightarrow \Epp^\infty_{0,q} = H_q(G,C),
\]
induced by the projections $(f\otimes_G c)\mapsto (f\otimes_G c)$.

\section{Lower triangle lemma}

For the proofs of the two main theorems we will need a technical lemma on convergence of a spectral sequence with a certain `zero triangle'.

\begin{Lemma}[Lower triangle lemma]
	Let $(E^r_{p,q})_{r\geq 0}$ be a spectral sequence converging to $H_{p+q}$. Assume that the second page $E^2$ has the following structure:
	\[
		E^2_{p,q} = \left\{\begin{array}{ccc}
			G_q & p=0,& 0\leq q \leq n \\
			0 & p>0, &p+q\leq n
		\end{array}\right..
	\]
	Then the filtration on $H$ is trivial, namely,
	\[
		\Gr_{p,q}H = \left\{\begin{array}{ccc}
			H_q & p=0,& 0\leq q \leq n \\
			0 & p>0, &p+q\leq n
		\end{array}\right.
	\]
	and consequently the edge homomorphism
	\[
		\pi_q: G_q = E^2_{0,q} \twoheadrightarrow E^\infty_{0,q} = H_q
	\]
	is surjective for $0\leq q\leq n$. It is an isomorphism for $0\leq q < n$. \label{l:triangle}
\end{Lemma}

\begin{figure}[hbt]
	\input{triangle} \caption{The conditions of the lower triangle lemma} \label{f:triangle}
\end{figure}

\begin{Proof}
	Since this type of edge homomorphism is always surjective, there is not much to prove. All groups which are already zero will stay zero on every page. Hence $H_q = \cF^0H_q = \Gr_{0,q} H = E^\infty_{0,q}$ for all $q\leq n$, the filtration is trivial.
	
	On the other hand, for all $E^k_{0,q}$ with $q<n$ all differential arrows begin and end in zero groups, hence $\pi_q: E^2_{0,q} \twoheadrightarrow E^\infty_{0,q}$ is an isomorphism for $q<n$.
\end{Proof}

\pagebreak
\section{The Lyndon/Hochschild-Serre spectral sequence}
The following spectral sequence is a tool to calculate homology of groups via their normal subgroups.\nopagebreak
\begin{BigTheorem}[Lyndon/Hochschild-Serre spectral sequence]
	For any normal subgroup $H$ of a group $G$ and for any $G$-module $M$, there is a convergent first-quadrant spectral sequence of homological type as follows:
	\[
	E^2_{p,q}=H_p\bigl(G/H, H_q(H,M)\bigr) \quad\Rightarrow\quad H_{p+q}(G,M).
	\]
	The action of $G/H$ on $H_q(H,M)$ is as in Corollary~\ref{c:ghaction}. It is induced by
	\[
		gH \cdot (f\otimes_H m) = (c_g(f)\otimes_H gm),
	\]
	where $c_g$ is the diagonal conjugation on $F_q(H)$.
	The edge maps 
	\begin{align*}
		H_p(G,M) &= E^\infty_{p,0} &\stackrel{\iota_p}{\hookrightarrow} &E^2_{p,0} = H_p(G/H, H_0(H,M))= &H_p(G/H, \Z \otimes_{\Z H} M) \\
		\intertext{are induced by the map}
		(f\otimes_G m) &&\mapsto &&( p_1(f) \otimes_{G/H} p_2(m) ),
	\end{align*}
	where
	\[
		p_1: G \twoheadrightarrow G/H,\qquad p_2: M \twoheadrightarrow \Z \otimes_{\Z H} M
	\]
	are the canonical projections. \label{t:hochschildserre}
\end{BigTheorem}
\nopagebreak
\begin{Proof}
	For a proof, see \cite[Theorem~6.8.2]{Wbl:IHA:94} or \cite[Theorem~VII.6.3]{Bro:CoG:82}.
\end{Proof}

\chapter{Homological isomorphisms via inclusions} \label{c:homisom}
In this chapter, we state a technical result by A.~Suslin, which will be used later on. It can be found in \cite[\S\space 1]{Sus:HoG:84}. To prove the result, we first have to show some intermediate results about the homology of multiplicative groups of fields.

\section{(Non-)Existence of certain field embeddings}

We will consider the following situation:

\begin{Problem} Given two fields $F$ and $F'$, when do field embeddings
	\[
		\varphi_1,\ldots,\varphi_n: F \hookrightarrow F'
	\]
	such that for any $x \in F^\star$ we have
	\[
		\prod_{i=1}^n \varphi_i(x) = 1
	\]
	exist? \label{p:fieldembeddings}
\end{Problem}
\begin{Lemma} If in the situation~\ref{p:fieldembeddings} the field $F$ is finite with $|F|=p^m$, then we necessarily have
	\[
		n\geq (p-1)\cdot m.
	\] \label{l:finitefieldemb}
\end{Lemma}
\begin{Proof}
The field generated by all subfields $\im\varphi_i$ is a finite subfield of $F'$, since it is the vector space generated by all products of image elements, this is a finite-dimensional vector space over a finite field.

Every finite field admits at most one subfield with $p^m$ elements (see for example \cite[Corollary~1 of Theorem~4.26]{Jac:BaI:85}). So we know that $\im\varphi_i =\im\varphi_1$ for all $i$ and $\im\varphi_i \cong F$. We can hence assume that $F'=F$ and that all $\varphi_i$ are field automorphisms of $F$.
	
For any finite field, we know that $\Aut(F)=\langle\sigma\rangle$ (see for example \cite[Theorem~4.26]{Jac:BaII:89}, where $\sigma:F\rightarrow F$, $x\mapsto x^p$ is the Frobenius automorphism. So there exist integers $t_i\in \Z$ such that $\varphi_i = \sigma^{t_i}$ for all $1\leq i\leq n$. 
	
	\noindent Then we know that
	\begin{eqnarray*}
		\prod_{i=1}^n x^{p^{t_i}}  &=& 1\quad\text{for all $x\in F^\star$, which means}\\
		\sum_{i=1}^n p^{t_i} &\equiv& 0 \mod (p^m-1).
	\end{eqnarray*}
	We denote by $n_j$ the number of $t_i$ such that $t_i \equiv j \mod m$. Then $\sum_{j=0}^{m-1} n_j = n$ and
	\begin{equation}
		\sum_{j=0}^{m-1} n_jp^j \equiv 0 \mod (p^m-1) \label{e:numbersublemma}.
	\end{equation}
	\subparagraph{Claim} For all non-negative integers $n_j$ which are not all zero and satisfy \eqref{e:numbersublemma}, we have 
		\[
			\sum_{j=0}^{m-1} n_j\geq (p-1)\cdot m.
		\]	
	\subparagraph{Proof}
		If $n_j \geq p$ for a certain $j$, by replacing $n_j$ with $n_j-p$ and $n_{j+1}$ by $n_{j+1}+1$ we get a new system of integers $n_j$ that satisfy equation \eqref{e:numbersublemma} while diminishing $\sum_j n_j$. Thus we may suppose that $0 \leq n_j <p$ for all $j$. Then 
		\[
			0 < \sum_{j=0}^{m-1} n_j\cdot p^j \leq (p-1)(1+p+\cdots+p^{m-1})=p^m-1.
		\]
		By \eqref{e:numbersublemma} we know that
		\[
			\sum_{j=0}^{m-1}n_j\cdot p^j = p^m-1.
		\]
		It follows that $n_j=p-1$ for all $j$ and hence
		\[
			\sum_{j=0}^{m-1}n_j=(p-1)\cdot m.
		\]
	So we know that $n\geq (p-1)\cdot m$. 
\end{Proof}

\begin{Lemma}
	If in the situation~\ref{p:fieldembeddings} the field $F$ is infinite, such field embeddings cannot exist. \label{l:inffieldemb}
\end{Lemma}\nopagebreak
\begin{Proof} We assume that these field embeddings exist and distinguish the following cases:

	If $F$ has characteristic zero, then for some $n\geq 1$ we have the equation
	\[
	2^n = \prod_{i=1}^n \varphi_i(2) = 1,
	\]
	which is impossible.
	
	If $F$ is an infinite algebraic extension of a finite field, it contains arbitrarily big finite subfields. We restrict the field embeddings to these subfields. We can apply Lemma~\ref{l:finitefieldemb}. Then $n\geq (p-1)m$, where $m$, the size of the subfield, grows arbitrarily large. This is impossible.

	Now, assume that $F$ is transcendental over its prime subfield $\F_p$, let $t\in F$ be transcendental over $\F_p$ and set $t_i := \varphi_i(t)$. Then every $t_i$ is transcendental over $\F_p$, the prime field of $F'$. For any non-zero polynomial $f\in \F_p[X]$, we have
	\[
		\prod_{i=1}^n f(t_i) = \prod_{i=1}^n \varphi_i(f(t)) = 1,
	\]
since $f(t)\neq 0$ for any polynomial $f\in\F_p[X]$ and the embeddings $\varphi_i$ map $\F_p$ identically to $\F_p$. We form the ideal generated by all these polynomials:
\[
	I(\F_p) := \langle \bigl(\prod_{i=1}^n f(X_i)\bigr) -1 \,|\, 0\neq f \in \F_p[X] \rangle \subseteq \F_p[X_1,\ldots,X_n].
\]
Now, to every field extension $E$ of $\F_p$, we assign an affine variety $V(E)\subseteq \A^n(E)$ defined by the ideal
\[
	I(E) := E\cdot I(\F_p) = \langle \bigl(\prod_{i=1}^n f(X_i)\bigr) -1 \,|\, 0\neq f \in \F_p[X] \rangle \subseteq E[X_1,\ldots,X_n].
\]	
We just proved that $V(F')$ is non-empty, it contains $(t_1,\ldots,t_n)$. This of course implies that $I(F')\neq F'[X_1,\ldots,X_n]$, hence $I(\F_p)\neq \F_p[X_1,\ldots,X_n]$.

Let now $K$ be the algebraic closure of $\F_p$, then we also know that 
\[
	I(K)\neq K[X_1,\ldots,X_n].
\]
Then by a special case of Hilbert's Nullstellensatz (see \cite[Theorem~7.17]{Jac:BaII:89}), $V(K)$ is also non-empty.

But every point in $(x_1,\ldots,x_n)\in V(K)$ has coordinates which are algebraic over $\F_p$. Hence we can choose a point $(x_1,\ldots,x_n)\in V(K)$ and a polynomial $f\in \F_p[X]$ such that $f(x_1)=0$. Then we get
	\[
		0 = \prod_{i=1}^n f(x_i) = 1,
	\]
	which is a contradiction. The field embeddings cannot exist.
\end{Proof}

\section{Homology of the multiplicative groups of fields}
In the following, for any $n\geq 1$ and any field $F$, we write $F^{\otimes n} = F\otimes\cdots\otimes F$, which is a commutative $F$-algebra with pointwise multiplication.
\begin{Corollary}
	If $F$ is an infinite field and $F^{\otimes n}$ is endowed with the diagonal action of $F^\star$, then $H_0(F^\star,F^{\otimes n})=0$. \label{c:h0fstarfn}
\end{Corollary}
\begin{Proof}
	As in the proof of Lemma~\ref{l:vanishinghomology}, we know that $H_0(F^\star, F^{\otimes n}) = F^{\otimes n}/ I$, where $I$ is the ideal generated as follows:
	\[
	I=\langle(f\otimes\cdots\otimes f) - (1\otimes\cdots\otimes 1)\,|\, f\in F^\star\rangle.
	\]
	If $I\neq F^{\otimes n}$, then $I$ it is contained in a maximal ideal $J$. Then the maps 
	\[
		F\rightarrow F^{\otimes n},\quad f\mapsto (1\otimes \cdots \otimes f \otimes \cdots \otimes 1),
	\]
	which embed $F$ into the tensor product induce field embeddings $F\rightarrow F^{\otimes n}/J$. These fulfil the condition of Lemma~\ref{l:inffieldemb}, which is impossible since $F$ is infinite.
\end{Proof}
\begin{Definition}
	We denote by $V_n$ the subalgebra of $F^{\otimes n}$ of all symmetric tensors (those which are invariant under the action of $S_n$). Then $F^{\otimes n}$ is a $V_n$-module via left multiplication in $F^{\otimes n}$.
\end{Definition}
\begin{Proposition}
	For any infinite field $F$, we have $H_0(F^\star,V_n)=0$. \label{p:homologyVn}
\end{Proposition}
\begin{Proof}
	Denote by $V_n^{n_1,\ldots,n_r}\subseteq V_n$ the subspace of $V_n$ generated by the tensors of the following form:
	\[
	\sum_{\sigma\in S_n/S_{n_1}\times\cdots\times S_{n_r}} \bigl(\underbrace{f_1\otimes\cdots\otimes f_1}_{n_1 \text{ times}}\otimes \cdots \otimes \underbrace{f_r\otimes\cdots\otimes f_r}_{n_r \text{ times}}\bigr)^\sigma,
	\]
	where the $n_i$ are positive integers with $\sum_{i=1}^r n_i = n$. These $V_n^{n_1,\ldots,n_r}$ are $F^\star$-in\-vari\-ant subspaces, hence $F^\star$-submodules of $V_n$ and $V_n = \sum V_n^{n_1,\ldots,n_r}$. We define an $F^\star$-in\-vari\-ant filtration on $V_n$ by setting:
	\[
		V_n^{(i)} = \sum_{\sum(n_j-1)\leq i} V_n^{n_1,\ldots,n_r}.
	\]
	Let $0\leq i\leq n$. For any $n_1,\ldots,n_r$ satisfying $\sum(n_j-1)= i$, the map
	\begin{eqnarray*}
		F^r &\rightarrow& V_n^{(i)}/V_n^{(i-1)}\\
		(f_1,\ldots,f_r)&\mapsto& \sum_{\sigma\in S_n/S_{n_1}\times\cdots\times S_{n_r}}\bigl(\underbrace{f_1\otimes\cdots\otimes f_1}_{n_1 \text{ times}}\otimes \cdots \otimes \underbrace{f_r\otimes\cdots\otimes f_r}_{n_r \text{ times}}\bigr)^\sigma \\ &&+ V_n^{(i-1)}
	\end{eqnarray*}
	is multilinear and hence induces a homomorphism of $F^\star$-modules
	\[
	F^{\otimes r} \rightarrow V_n^{(i)}/V_n^{(i-1)}.
	\]
	Forming direct sums over all possible values of the $n_j$, we get an epimorphism
	\[
	\bigoplus_{\genfrac{}{}{0pt}{}{n_1,\ldots,n_r}{\sum(n_j-1)=i}} F^{\otimes r} \twoheadrightarrow V_n^{(i)}/V_n^{(i-1)}.
	\]
	This is automatically a split epimorphism, since these are all vector spaces. We apply the homology functor. Since this is a split epimorphism, we get an induced epimorphism on homology modules:
	\[
		H_0(F^\star,\bigoplus_{\genfrac{}{}{0pt}{}{n_1,\ldots,n_r}{\sum(n_j-1)=i}} F^{\otimes r}) \twoheadrightarrow H_0\bigl(F^\star,V_n^{(i)}/V_n^{(i-1)}\bigr).
	\]
	By the additivity of the homology functor in the coefficient module and by Corollary~\ref{c:h0fstarfn}, we now know that
	\[
		H_0\bigl(F^\star,V_n^{(i)}/V_n^{(i-1)}\bigr)=0.
	\]
	We write down the long exact homology sequence associated to the short exact sequence
	\[
		0 \rightarrow F(F^\star)\otimes V_n^{(i-1)} \rightarrow F(F^\star)\otimes V_n^{(i)} \rightarrow F(F^\star)\otimes V_n^{(i)}/V_n^{(i-1)} \rightarrow 0
	\]
	of chain complexes of vector spaces (where $F(F^\star)$ is the standard resolution of $F^\star$):
	\[
		\cdots\rightarrow H_0(F^\star,V_n^{(i-1)}) \rightarrow H_0(F^\star,V_n^{(i)}) \rightarrow H_0(F^\star,V_n^{(i)}/V_n^{(i-1)}) \rightarrow 0
	\]
	By a simple induction on $i$, starting with $i=0$, we deduce that $H_0(F^\star,V_n^{(i)})=0$ for all $0\leq i \leq n$, from which it follows that $H_0(F^\star,V_n)=0$.
\end{Proof}

\begin{Lemma}
	Let $n_1,\ldots,n_r$ be positive integers and for any $1\leq i\leq r$ let $T^{n_i}$ be either the external or symmetric power of an infinite field $F$ over $\Z$ or the prime subfield, that is
	\[
	T^{n_i}\in\bigl\{\Lambda^{n_i}(F),S^{n_i}(F)\bigr\}.
	\]
	Put $n=\sum n_i$ and consider the projection
	\[
		F^{\otimes n} \rightarrow T^{n_1}\otimes \cdots \otimes T^{n_r}.
	\]
	Then the kernel $K$ of this projection is a $V_n$-submodule of $F^{\otimes n}$ and hence
	\[
	T^{n_1}\otimes \cdots \otimes T^{n_r}\cong F^{\otimes n}/K
	\]
	has the natural structure of a $V_n$-module. \label{l:hompowers}
\end{Lemma}

\begin{Proof}
	Let $K_i$ denote the kernel of the map $F^{\otimes n_i}\rightarrow T^{n_i}$. Then 
	\[
		K=\sum_{i=1}^r F^{\otimes n_1} \otimes\cdots\otimes K_i \otimes\cdots\otimes F^{\otimes n_r}
	\]
	and it suffices to check that each summand is a $V_n$-submodule of $F^{\otimes n}$. We may assume that $i=1$. We may further assume that $n_1\geq 2$, since $K_1=0$ for $n_1=1$.

\begin{description}
	\item[Case 1: $\mathbf{T^{n_1}=S^{n_1}(F)}$ is the symmetric power] Then $K_1\otimes\cdots\otimes F^{\otimes n_r}$ is generated by tensors of the form $v^{(i,j)}-v$, where $1\leq i<j\leq n_1$. Then, for any $w\in V_n$, we calculate
	\[
		w\cdot(v^{(i,j)}-v) = (wv)^{(i,j)} - (wv) \in K_1\otimes\cdots\otimes F^{\otimes n_r},
	\]
	hence $K_1\otimes\cdots\otimes F^{\otimes n_r}$ is closed under scalar multiplication with $V_n$.
\item[Case 2: $\mathbf{T^{n_1}=\Lambda^{n_1}(F)}$ is the exterior power] In this case, $K_1\otimes\cdots\otimes F^{\otimes n_r}$ is generated by tensors $v$, such that $v=v^{(i,j)}$ for some $1\leq i<j\leq n_1$. Then, for any $w\in V_n$, we see that
	\[
		(wv)^{(i,j)}=w^{(i,j)}v^{(i,j)}=wv,
	\]
	and hence $wv\in K_1\otimes\cdots\otimes F^{\otimes n_r}$, which is thus also a submodule.
	\end{description}
\end{Proof}

\begin{Corollary}
	Under the conditions of Lemma~\ref{l:hompowers}, for any $i\geq 0$ we have 
	\[
	H_i\bigl(F^\star,T^{n_1}\otimes \cdots \otimes T^{n_r}\bigr)=0.
	\] \label{c:hompowerszero}
\end{Corollary}
\begin{Proof}
	Set $A=V_n$ and $M=T^{n_1}\otimes \cdots \otimes T^{n_r}$, which is an $A$-module by Lemma~\ref{l:hompowers}. By Proposition~\ref{p:homologyVn} we know that $H_0(F^\star,A)=0$. So we can apply Lemma~\ref{l:vanishinghomology} to see that $H_i(F^\star,M)=0$ for any $i\geq 0$. The group homomorphism $\varphi$ for Lemma~\ref{l:vanishinghomology} is the diagonal inclusion $F^\star\hookrightarrow (F^{\otimes n})^\star$.
\end{Proof}

\begin{Proposition}
	Suppose that $F$ is an infinite field of characteristic $p>0$ and $V$ is a vector space over $F$. Then for any $i\in\Z$ and any $j>0$ we have 
	\[
	H_i\bigl(F^\star, H_j(V, \F_p)\bigr)=0.
	\]
	Here, the action of $F^\star$ on $H_j(V,\F_p)$ is induced by the diagonal left multiplication action on $F(V)$. \label{p:homvectfp}
\end{Proposition}
\begin{Proof}
	By the Corollaries~\ref{c:homvectp} and \ref{c:homvect2} we know that for any $\F_p$-vector space $V$ the graded $\F_p$-vector space $H_\stern(V,\F_p)$ is isomorphic to $\Lambda(V)\otimes S(V)$ if $p\neq 2$ and to $S(V)$ if $p=2$.

Hence $H_j(V,F_p)$ is isomorphic as an $F^\star$-module to a direct sum of modules considered in Lemma~\ref{l:hompowers} and the statement follows from Corollary~\ref{c:hompowerszero}.
\end{Proof}

\begin{Theorem}
	Let $F$ be any infinite field, $V$ an $F$-vector space and $k$ any prime field. Then for any $i\in\Z$ and any $j>0$ we have
	\[
		H_i\bigl(F^\star, H_j(V, k)\bigr)=0.
	\] \label{t:vanish}
\end{Theorem}
\begin{Proof} We have to consider the following cases.
	\begin{description}
		\item[Case 1] If $\chr k \neq \chr F$, then $F_j(V)\otimes_V k=0$, hence $H_j(V,k)=0$, hence
			\[
				H_i(F^\star, H_j(V,k))=0.
			\]
		\item[Case 2] If $\chr k = \chr F > 0$, then the statement follows from Proposition~\ref{p:homvectfp}.
		\item[Case 3] If $\chr k = \chr F = 0$, then by Corollary~\ref{c:homvectq} we know that $H_j(V,k)=\Lambda^j(V)$. For any positive integer $n$ the action of $n\in F^\star$ on $H_j(V,k)$ coincides with multiplication by $n^j$, since it is induced by diagonal multiplication on $F_j(V)$. But $F^\star$ is abelian, by Proposition~\ref{p:conj} we know that $H_i(F^\star, H_j(V,k))$ is annihilated by $n^j-1$ and hence must be zero, since it is a $\Q$-vector space.
	\end{description}
\end{Proof}

\section{A homological isomorphism theorem}

\begin{Theorem}[Suslin]
	Let $F$ be any skew-field with infinite centre $K$, let $m,n\geq 1$ and let $G_1\leq \Gl_n F$, $G_2\leq \Gl_m F$ be subgroups such that at least one of then contains the group of $K^\star$-multiples of the identity matrix. Let $M\subseteq M_{n,m}(F)$ be a vector subspace such that $G_1 M = M = MG_2$. Then the natural embedding
	\[
		i: \begin{pmatrix} G_1 & 0 \\ 0 & G_2 \end{pmatrix} \hookrightarrow \begin{pmatrix} G_1 & M \\ 0 & G_2 \end{pmatrix}
	\]
	induces an isomorphism on the corresponding homology groups. \label{t:suslin}
\end{Theorem}
\begin{Proof}
	We abbreviate $G=\GMG$ and we know that $\GoG\cong G_1\times G_2$ and $\oMo$ is isomorphic to the additive group of $M$. Furthermore it is clear that
	\[
		G/\oMo \cong G_1\times G_2.
	\]
	Let $k$ be any prime field. We consider the Lyndon/Hochschild-Serre spectral sequence $E$ (Theorem~\ref{t:hochschildserre}) associated to the normal subgroup $\oMo$ in $G$. The spectral sequence has the form
	\[
		E^2_{p,q} = H_p\Bigl(\GoG, H_q\bigl(\oMo,k\bigr)\Bigr) \Rightarrow H_{p+q}(G,k).
	\]
	Here, the action of $\GoG$ on $H_q(\oMo,k)$ is induced by diagonal conjugation in the standard resolution, the action on $k$ is trivial. We calculate
	\[
		\begin{pmatrix} g_1 & 0 \\ 0 & g_2\end{pmatrix}\begin{pmatrix} I_n & m \\ 0 & I_m\end{pmatrix}\begin{pmatrix} g_1^{-1} & 0 \\ 0 & g_2^{-1}\end{pmatrix} = \begin{pmatrix} I_n & g_1mg_2^{-1} \\ 0 & I_m\end{pmatrix}.
	\]
	We apply the isomorphisms $G_1\times G_2\cong\GoG$ and $\oMo\cong M$ (additively). Then we can write down the spectral sequence more efficiently as
	\[
		E^2_{p,q} \cong H_p(G_1\times G_2, H_q(M,k)) \Rightarrow H_{p+q}(G,k).	
	\]
We know that the action of $G_1\times G_2$ on $H_q(M,k)$ is induced by
	\[
		(g_1,g_2)\cdot ( (m_0,\ldots,m_q)\otimes 1) = ( (g_1m_0g_2^{-1},\ldots,g_1m_qg_2^{-1})\otimes 1)
	\]
and linear extension.

Our first goal will be to show that $E^2_{p,q} = 0$ for $q\neq 0$. We can assume that, without loss of generality, $G_1$ contains the subgroup $K^\star I_n\cong K^\star$, where $I_n$ is the identity matrix. This subgroup is central, hence normal. We then have 
	\[
		L := K^\star I_n \times \{I_m\} \trianglelefteq G_1 \times G_2
	\]
	and for any $q$ we get another Lyndon/Hochschild-Serre spectral sequence $\Ep$ as follows:
	\[
	\Ep^2_{i,j} = H_i((G_1\times G_2)/L,H_j(L,H_q(M,k))) \Rightarrow H_{i+j}(G_1\times G_2, H_q(M,k)) = E^2_{i+j,q}
	\]
	Note that, by what we have said above, the action of $L$ on $H_q(M,k)$ is induced by
	\[
		(l\cdot I_n,I_m)\cdot ( (m_0,\ldots,m_q)\otimes 1) = ( (lm_0,\ldots,lm_q)\otimes 1)
	\]
	and linear extension. So the action of $L$ on $H_q(M,k)$ coincides with the action of $K^\star$ on $H_q(M,k)$. By Theorem~\ref{t:vanish} we then know that $H_j(L,H_q(M,k))=0$ for $q\neq 0$. This in turn means that $\Ep^2_{i,j}=0$ if $q\neq 0$, which implies that the abutment vanishes as well.

	Hence the first spectral sequence $E$ collapses at the second page, since $E^2_{p,q}=0$ for $q\neq 0$ and all further differentials begin and end in trivial groups. The edge map
	\[
	H_p(G,k) = E^\infty_{p,0}\stackrel{\iota_p}{\hookrightarrow} E^2_{p,0} = H_p\Bigl(\GoG,\underbrace{H_0\bigr(\oMo,k\bigl)}_{\cong k}\Bigr)
	\]
	is an injective homomorphism induced by the projection
	\[p_1: G \rightarrow G/\oMo\cong\GoG,\] which is split by the inclusion $i:\GoG \rightarrow G$ in the other direction, hence the edge map $\iota_p$ is also surjective.

	Then the inclusion $i: \GoG \hookrightarrow G$ induces a homology isomorphism
	\[
	H_p(\GoG,k) \rightarrow H_p(G,k)
	\]
	for any prime field $k$. By a classical theorem in homological algebra (see for example \cite[Corollary~3A.7.~(b)]{Hat:ATp:02}), $i$ must then induce a homomorphism on homology with coefficients in $\Z$, which proves the theorem.
\end{Proof}

\chapter{Two homological stability theorems} \label{c:homstabthms}
The main topics of this thesis are two theorems on homological stability --- one for standard unitary groups over $\R$, $\C$ and $\bH$, and one for general linear groups over skew-fields with infinite centres. These theorems can both be proven with analogous methods. In this chapter, we will try to explain this method and examine similar and different steps for both cases.

This proof method and the two resulting theorems are due to Chih-Han Sah (see \cite[Theorems~1.1 and B.1]{Sah:HcL:86}).

Another proof of the general linear group stability (over fields) can be found in \cite[Theorem~3.4]{Sus:HoG:84}. Homological stability results (with weaker stability ranges) are also known for general linear and unitary groups over more general rings, see for example \cite{Wag:SfH:76}, \cite{vdK:HSL:80} and \cite{MaB:HSU:02}.

\section{The setting}

We will be dealing with many different inclusions in this chapter, hence we introduce a convenient notation. 

\paragraph{Notation} If we have an inclusion $A\hookrightarrow B$, we will denote the inclusion map by \[i^A_B: A \hookrightarrow B.\]

\subsection{The theorems}

We will consider the following situations:

\ruleboxstart
\minisec{The unitary case} Let $\F$ be either $\R$, $\C$ or $\bH$ and $\star$ be the standard involution of $\F$. For $n\geq 1$, we endow $\F^n$ with the standard hermitian positive definite inner product as follows:
\[
	\langle u, v \rangle = \sum_i u_i^\star v_i.
\]
Let $G_n=U_n\F$ be the associated standard unitary group acting on $\F^n$. We pick any bases of $\F^n$ and $\F^{n+1}$. Then there is an inclusion
	\begin{align*}
	i^{U_n\F}_{U_{n+1}\F}: U_n \F &\rightarrow U_{n+1}\F \\
	A &\mapsto \bigl(\begin{smallmatrix} A & 0 \\ 0 & 1\end{smallmatrix}\bigr).
	\end{align*}
		
\begin{BigTheorem}[Homological stability, Sah]
The induced homomorphism
	\[
	(i^{U_n\F}_{U_{n+1}\F})_*: H_q(U_n \F) \rightarrow H_q(U_{n+1}\F)
	\]
	is surjective for $q\leq n$. It is an isomorphism for $q<n$. \label{t:unstability}
\end{BigTheorem}

\minisec{The general linear case} For $n\geq 1$, let $G_n=\Gl_n\F$ be the $n$-dimensional general linear group over $\F$, where $\F$ is any skew-field with infinite centre. We choose bases of $\F^n$ and $\F^{n+1}$. Then there is an inclusion
	\begin{align*}
	i^{\Gl_n\F}_{\Gl_{n+1}\F}: \Gl_n\F &\rightarrow \Gl_{n+1}\F.\\
	A &\mapsto \bigl(\begin{smallmatrix} A & 0 \\ 0 & 1\end{smallmatrix}\bigr).
	\end{align*}
	
\begin{BigTheorem}[Homological stability, Sah]
The induced homomorphism
	\[
	(i^{\Gl_n\F}_{\Gl_{n+1}\F})_*: H_q(\Gl_n \F) \rightarrow H_q(\Gl_{n+1}\F)
	\]
	is surjective for $q\leq n$. It is an isomorphism for $q<n$. \label{t:glnstability}\ruleboxend
\end{BigTheorem}

\paragraph{Basis independence} Note that we have picked the bases arbitrarily. But if we choose different bases, the resulting inclusion is conjugate via $U_{n+1}\F$ or $\Gl_{n+1}\F$, respectively, to any other inclusion with respect to another basis. By Proposition~\ref{p:conj}, these conjugations induce the identity on homology modules, hence we can choose any fixed basis without changing the induced map on homology.

\paragraph{Ambiguity of notation}
We want to tackle the proofs of these two theorems in parallel. Using a slight abuse of notation, we denote both sequences of groups by $G_n$. This notation will only be used in steps where both cases are completely analogous, and will (hopefully) not lead to confusion. Whenever we distinguish both cases, the corresponding section is enclosed between horizontal lines, as above.

\subsection{Prerequisite: Witt's theorem}

For the unitary case we will frequently use a consequence of a classical theorem by Witt. In the general linear case the corresponding result is almost trivial.

\ruleboxstart
\minisec{The unitary case}
We state a classical theorem by Witt:
\begin{BigTheorem}[Witt's extension theorem]
Let $V$ be a finite dimensional vector space endowed with a non-degenerate hermitian form. Any isometry of subspaces $U_1\rightarrow U_2$ an be extended to a unitary map $V\rightarrow V$.
\end{BigTheorem}
\begin{Proof}
See for example \cite[Theorem~7.4]{Tay:GCG:92}. See \cite[Theorem~6.2.12]{HoM:CGK:89} for a very general version of the theorem.
\end{Proof}
We apply this theorem to our situation and derive a special case.
\begin{Corollary} For $\F\in\{\R,\C,\bH\}$ as above, we pick vectors
\[
v_1,\ldots,v_k,w_1,\ldots,w_k\in \F^n
\]
with 
\[
\langle v_i, v_j\rangle=\langle w_i, w_j \rangle \quad\text{for all}\quad 1\leq i,j \leq k.
\]
We define
\[
V=\spn(v_1,\ldots,v_k),\quad W=\spn(w_1,\ldots,w_k).
\]
Then there exists a unitary map $\sigma\in U_n \F$ fixing $(V+W)^\perp$ with
\[
\sigma(v_i)=w_i \quad\text{for all}\quad 1\leq i\leq k.
\] \label{c:wittun}
\end{Corollary}
\begin{Proof}
By hypothesis, there is a unitary map $\sigma':V\rightarrow W$ with $\sigma'(v_i)=w_i$. It can be extended by the identity map on $(V+W)^\perp$ to a unitary map 
\[
	V\oplus (V+W)^\perp\rightarrow W\oplus(V+W)^\perp.
\]
We apply Witt's extension theorem to get a unitary map $\sigma:\F^n\rightarrow \F^n$ that has the required properties.
\end{Proof}\pagebreak
\minisec{The general linear case}
This is an analogue of the above corollary for the general linear case.\nopagebreak
\begin{Proposition}
	Given any $\F$-vector space $E$, for two $k$-tuples of linearly independent vectors
	\[
	v_1,\ldots,v_k\in E,\qquad w_1,\ldots,w_k\in E
	\]
	we set
	\[
		V=\spn(v_1,\ldots,v_k),\quad W=\spn(w_1,\ldots,w_k).
	\]
	Let $U$ be a complement of $V+W$, i.e.~$(V+W)\oplus U = E$. Then there exists an automorphism $\sigma\in \Gl (E)$ fixing $U$ with 
	\[
		\sigma(v_i)=w_i \quad\text{for all}\quad 1\leq i\leq k.
	\] \label{p:wittgln}
\end{Proposition}
\begin{Proof}
	Choose a basis $\{e_i\}$ of $V\cap W$, extend it to bases $\{e_i\}\cup\{v_j\}$ of $V$ and $\{e_i\}\cup\{w_l\}$ of $W$. Then $E=\langle e_i\rangle\oplus \langle v_j \rangle\oplus\langle w_l\rangle \oplus U$ and it suffices to write down an appropriate matrix in the finite-dimensional vector space
	\[
		\langle e_i\rangle\oplus \langle v_j \rangle\oplus\langle w_l\rangle=V+W.
	\]
\end{Proof}\vspace{-1em}\ruleboxend

\section{A spectral sequence}
The rest of this chapter will cover the proofs of the two theorems. The main idea in the proofs will be to use Corollary~\ref{c:acyclicC}, which gives us a spectral sequence whose abutment is isomorphic to the homology groups of $G_{n+1}$, given an acyclic $G_{n+1}$-complex. We will then be able to work in this spectral sequence.

We will start with many general observations. At the end the proofs will be by induction. For now, we fix an $n\geq 0$ and we further abbreviate $G=G_{n+1}$.

\subsection{Construction of a spectral sequence}
We already know a chain complex which is always acyclic: the ordered simplicial chain complex $C(X)$, where $X$ is a set with a left $G$-action, see Lemma~\ref{l:acycl}.

\ruleboxstart
\minisec{The unitary case}
We denote by $S^n$ the set of all vectors in $\F^{n+1}$ with norm $1$. We simply take $X$=$S^n$, on which $U_{n+1}\F$ acts transitively from the left. 

\minisec{The general linear case}
This is even simpler, we take $X=\F^{n+1}\backslash \{0\}$, the non-zero vectors of $\F^{n+1}$. Here we also have a transitive left $\Gl_{n+1}\F$-action.
\ruleboxend

\paragraph{}In both cases, we take $C := C(X)$, this is an acyclic chain complex of $G$-modules. Note that, again, the notation is ambiguous, we denote both chain complexes by $C$. Again, it will either be clear in which context we are or the argument will be completely analogous in both cases.

On these chain complexes, we will need a concept of a dimension. Since in both cases, $X$ is a subset of $\F^{n+1}$, we can give a definition in parallel.

\begin{Definition}
	For any element $c=\sum_{j\in J}(x^j_0,\ldots,x^j_p) \in C_p(X)$, we define the \emph{dimension of $c$} to be
	\[
		|c| := \max_{j\in J} \dim\spn(x^j_0,\ldots x^j_p).
	\]
	Note that this induces a natural filtration on $C(X)$, which is compatible with the action of $G$. We denote this filtration by
	\[
		\cF^l C_p = \{ c\in C_p : |c|\leq l+1\}.
	\]
\end{Definition}

\paragraph{The spectral sequence} Since we have constructed an acyclic chain complex of $G$-modules, by Corollary~\ref{c:acyclicC} there is a spectral sequence
\[
	E^1_{p,q} = H_q(G,C_p) \quad\Rightarrow\quad H_{p+q}(G,C)\cong H_{p+q}(G),
\]
the differential $d^1$ being induced by $\id_{F(G)} \otimes_G d_C$. We have additional results concerning the isomorphisms and edge homomorphisms, but this will only be required later on.

We examine this spectral sequence more closely in a sequence of lemmas. We begin with a discussion of stabilizer subgroups.

\subsection{Homology of stabilizer subgroups}

The structure of the stabilizer subgroups of the $G$-action on $C_p$ of course depends on the two cases.

\ruleboxstart
\minisec{The unitary case}
We know that any $p$-simplex $c=(v_0,\ldots,v_p)$ with $|c|=k$ spans a $k$-dimensional subspace of $\F^{n+1}$. Hence its stabilizer subgroup $G_c\leq G$ is conjugate to
\[
	G_c \cong_\varphi \Bigl\{ \begin{pmatrix} I_k & 0 \\ 0 & g \end{pmatrix} : g\in U_{(n+1)-k}\F \Bigr\},
\]
where $I_k$ is the identity matrix of size $k$ and $\varphi$ is the conjugation automorphism. This subgroup is isomorphic to $U_{(n+1)-k}\F$ via the canonical inclusion. We consider the homomorphism chain
\[
	U_{(n+1)-k}\F \hookrightarrow \begin{pmatrix} I_k & 0 \\ 0 & U_{(n+1)-k}\F \end{pmatrix} \stackrel{\varphi}{\rightarrow} G_c.
\]
By Proposition~\ref{p:conj}, conjugation induces the identity on homology groups, the inclusion is an isomorphism. It is induced by the inclusion into the stabilizer subgroup. Hence the composite map induces an isomorphism on the homology groups
\[
	H_q(G_c) = (i^{U_{(n+1-|c|)}\F}_{G_c})_\stern H_q(U_{(n+1)-|c|}\F).
\]

\minisec{The general linear case}
For general linear groups, stabilizers are a bit more complicated. Pick a $p$-simplex $c=(v_0,\ldots,v_p)$ with $|c|=k$, which then also spans a $k$-dimensional subspace of $\F^{n+1}$. This time, for an appropriate choice of a basis, its stabilizer subgroup $G_c\leq G$ is
\[
	G_c = \Bigl\{ \begin{pmatrix} I_k & m \\ 0 & g \end{pmatrix} : m\in M_{k,(n+1-k)}\F, g\in\Gl_{(n+1)-k}\F \Bigr\},
\]
where $I_k$ is the identity matrix of size $k$. So $G_c$ is conjugate to the subgroup
\[
	G_c \cong_\varphi \begin{pmatrix} I_k & M_{k,(n+1-k)}\F \\ 0 & \Gl_{(n+1)-k} \F \end{pmatrix},
\]
where $\varphi$ is the appropriate conjugation automorphism. We consider the following maps:
\[
\Gl_{(n+1)-k}\F \stackrel{\cong}{\rightarrow} \begin{pmatrix} I_k & 0 \\ 0 & \Gl_{(n+1)-k}\F \end{pmatrix} \hookrightarrow \begin{pmatrix} I_k & M_{k,(n+1-k)}\F \\ 0 & \Gl_{(n+1)-k}\F \end{pmatrix} \stackrel{\varphi}{\rightarrow} G_c
\]
By Theorem~\ref{t:suslin}, the inclusion induces an isomorphism on homology. By Proposition~\ref{p:conj}, the conjugation $\varphi$ induces the identity map. In addition, note that the composition of all these maps is just the inclusion $i^{\Gl_{(n+1-|c|)}\F}_{G_c}$. So, transporting this result to homology, we have
\[
	H_q(G_c) = (i^{\Gl_{(n+1-|c|)}\F}_{G_c})_\stern H_q(\Gl_{(n+1)-|c|}\F).
\]
\vspace{-1em}\ruleboxend
\vspace{0.5em}
\begin{Lemma}
	In both cases we are left with the isomorphisms
\[
	H_q(G_c) = (i^{G_{(n+1-|c|)}}_{G_c})_\stern H_q(G_{(n+1)-|c|}).
\] \label{l:smallergn}
\end{Lemma}

\subsection{Appearance of the smaller groups}

We will try to get a grasp on the structure of $E^1$. To that end, we choose a set of representative simplices of $G$-orbits of $X^{p+1}$, denoted by $B_p$. Then, of course, the elements of $B_p$ generate $C_p$ as a $G$-module. We then have a decomposition:
	\[
	C_p = \bigoplus_{c\in B_p} (\Z G)c \cong \bigoplus_{c\in B_p} (\Z \otimes_{G_c} G).
	\]
\begin{Lemma}
	There is a homology isomorphism
	\begin{align*}
		\bigoplus_{c\in B_p} H_q(G_{n+1-|c|}) &\rightarrow H_q(G,C_p)=E^1_{p,q}\\
		\intertext{induced by}
		\sum_{c\in B_p} f_c \otimes_{G_{n+1-|c|}} z_c &\mapsto \sum_{c\in B_p} (i^{G_{n+1-|c|}}_G)f_c \otimes_G  z_c c.
	\end{align*} \label{l:specseqstruct}
\end{Lemma}
\paragraph{Remark} Most of the time, it will be enough to know that there is an isomorphism. The explicit description of the map inducing this isomorphism will be needed at the end of the proof, though, where we want to see that a particular homomorphism is induced by an inclusion $G_n\hookrightarrow G$.

\begin{Proof} We will write down a chain of maps from $\bigoplus_{c\in B_p} F_q(G_{n+1-|c|},\Z)$ to $F_q(G,C_p)$ in which every map induces a homology isomorphism, thereby proving the lemma. We begin with the isomorphism from Lemma~\ref{l:smallergn}:
	\begin{align*}
		\bigoplus_{c\in B_p} H_q(G_{n+1-|c|}) &\rightarrow \bigoplus_{c\in B_p} H_q(G_c)\\
		\intertext{induced by}
		\sum_{c\in B_p} f_c \otimes z_c &\mapsto \sum_{c\in B_p} (i^{G_{n+1-|c|}}_{G_c}) f_c \otimes_{G_c} z_c.
	\end{align*}
Next we apply Shapiro's Lemma (Lemma~\ref{l:shapiro}): The map
	\begin{align*}
		\cdots\rightarrow\bigoplus_{c\in B_p} H_q(G_c) &\rightarrow \bigoplus_{c\in B_p} H_q(G,\Z\otimes_{G_c} \Z G) \\
		\intertext{is an isomorphism induced by}
		\sum_{c\in B_p} (i^{G_{n+1-|c|}}_{G_c}) f_c \otimes_{G_c} z_c &\mapsto \sum_{c\in B_p} (i^{G_{n+1-|c|}}_G) f_c \otimes_G (z_c \otimes_{G_c} 1_G).
	\end{align*}
Since we know that $\Z\otimes_{G_c} \Z G$ is isomorphic to the $G$-orbit of $c$, hence to the $G$-module $\Z G c$ generated by $c$, we get
	\begin{align*}
		\cdots\rightarrow\bigoplus_{c\in B_p} H_q(G,\Z\otimes_{G_c} \Z G) &\rightarrow \bigoplus_{c\in B_p} H_q(G,\Z Gc)\\
		\intertext{induced by}
		\sum_{c\in B_p} (i^{G_{n+1-|c|}}_G) f_c \otimes_G (z_c \otimes_{G_c} 1_G) &\mapsto \sum_{c\in B_p} (i^{G_{n+1-|c|}}_G) f_c \otimes_G z_c1_Gc.
	\end{align*}
Finally we apply the additivity of the homology functor and get an isomorphism
	\begin{align*}
		\cdots\rightarrow\bigoplus_{c\in B_p} H_q(G,\Z Gc) &\rightarrow H_q(G,\bigoplus_{c\in B_p} \Z Gc)=H_q(G,C_p)\\
		\intertext{induced by}
		\sum_{c\in B_p} (i^{G_{n+1-|c|}}_G) f_c \otimes_G z_c1_Gc &\mapsto \sum_{c\in B_p} (i^{G_{n+1-|c|}}_G) f_c \otimes_G z_cc.
	\end{align*}
By concatenating all these maps, we get a proof of the lemma.
\end{Proof}
\paragraph{Leftmost column}In the special case of $p=0$, we know that $C_0 = \Z X$, on which $\Z G$ acts transitively. So $B_0=\{c\}$ consists of only one element in $X$. We get a homology isomorphism
	\begin{align*}
		H_q(G_n) &\rightarrow H_q(G,C_0)=E^1_{0,q}\\
		\intertext{induced by}
		f \otimes_{G_n} z &\mapsto (i^{G_n}_G)f \otimes_G z c.
	\end{align*}

This means that $E^1_{0,q}\cong H_q(G_n)$. The isomorphism is `almost' induced by the inclusion $G_n\hookrightarrow G$. We will see later on how this can be used to prove the theorem.

\paragraph{Bottom row} In the other special case of $q=0$, we know that $H_0(G_i)\cong \Z$ for any $i$. Hence we have an isomorphism
\begin{align*}
		\bigoplus_{c\in B_p} \Z &\rightarrow H_0(G,C_p)=E^1_{p,0}\\
		\intertext{induced by}
		\sum_{c\in B_p} z_c &\mapsto \sum_{c\in B_p} (1_G) \otimes_G  z_c c.
\end{align*}
We will denote an element $(1_G)\otimes_G 1\cdot c$ in $H_0(G,C_p)$ by $[c]$. The above isomorphism can then be written as
\begin{align}
		\Z B_p &\rightarrow H_0(G,C_p)=E^1_{p,0}\label{e:bottomrow}\\
		\intertext{induced by}
		\sum_{c\in B_p} z_c &\mapsto \sum_{c\in B_p} z_c [c].\nonumber
\end{align}
Hence we know that $E^1_{p,0}\cong \Z B_p$, the free abelian group over $G$-orbits of $X^{p+1}$.

The differential $d^1$ induces a differential $d^1$ on this {`}row' chain complex, we have $d^1[c] = [d_C c]$. Note that we have $[c]=[d]$ if and only if $c$ and $d$ are in the same $G$-orbit. In addition, for any simplex $c$, the dimension $|[c]|:=|c|$ is well defined, since $G$ does not change the dimension of simplices.

\subsection{Persistence of the smaller groups}

Since we want to apply the triangle lemma, we need to know that the leftmost column in $E^1$, containing modules isomorphic to $H(G_n)$, survives to the second page $E^2$. So we want to show
\begin{Lemma}
	$d^1_{1,q}: E^1_{1,q}\rightarrow E^1_{0,q}$ is the zero homomorphism, this implies
	\[
		E^2_{0,q}=E^1_{0,q} = H_q(G,C_0) \cong H_q(G_n).
	\]
\end{Lemma}

\begin{Proof}
	For this proof, we denote the boundaries in $(F(G)\otimes_G C_i)_q$ by $B_q(G,C_i)$. By Lemma~\ref{l:specseqstruct} we have an isomorphism
	\begin{align*}
		\bigoplus_{(v,w) \in B_1} H_q(G_{n-1}) &\rightarrow H_q(G,C_1)=E^1_{1,q}\\
		\intertext{induced by}
		\sum_{(v,w) \in B_p} f_{(v,w)} \otimes_{G_{n-1}} z_{(v,w)} &\mapsto \sum_{(v,w)\in B_p} (i^{G_{n-1}}_G)f_{(v,w)} \otimes_G  z_{(v,w)} (v,w),
	\end{align*}
	so it is enough to show that the differential maps all elements of the form
\[
	(i^{G_{n-1}}_G) f \otimes_G (v,w) + B_q(G,C_1)
\]
	onto the boundaries $B_q(G,C_0)$. 

	Pick any pair of elements $(v,w)\in X^2$. If $v$ and $w$ are linearly independent, we can choose a (possibly orthonormal) basis $(e_i)_{0\leq i \leq n}$ of $\F^{n+1}$, such that $v,w \in \spn(e_0,e_1)$. By Corollary~\ref{c:wittun} or, respectively, Proposition~\ref{p:wittgln}, there is an element $\sigma\in G$, that maps $v$ to $w$, leaving $\spn(e_2,\ldots e_n)$ fixed.

	If $v$ and $w$ are collinear, pick the basis such that $v=e_0$ and $\sigma\in G$ such that $\sigma(v)=w$ and the map $\sigma$ leaves $\spn(e_1,\ldots e_n)$ fixed.
	
	Remember that the element $f$ is a $(q+1)$-tuple of matrices in $G_{n-1}$. If we pick the inclusion into $G$ to be onto $\spn(e_2,\ldots e_n)$, which we can do by Proposition~\ref{p:conj}, then in both cases $\sigma$ commutes with every matrix in the tuple $(i^{G_{n-1}}_G) f$.

	We calculate the differential
\begin{equation*}
	\begin{split}
	d^1\bigl((i^{G_{n-1}}_G) f \otimes_G (v,w)&+B_q(G,C_1)\bigr)= (i^{G_{n-1}}_G) f \otimes_G d_C(v,w) +B_q(G,C_0)\\
	&=(i^{G_{n-1}}_G) f \otimes_G (w) - (i^{G_{n-1}}_G) f\otimes_G (v) +B_q(G,C_0)\\
	&=(i^{G_{n-1}}_G) f \otimes_G (w) - \sigma (i^{G_{n-1}}_G) f \sigma^{-1} \otimes_G (\sigma(v)) + B_q(G,C_0)\\
	&=0+B_q(G,C_0),
	\end{split}
\end{equation*}
where in the third line we conjugate with $\sigma$ and then apply Proposition~\ref{p:conj}, which implies that this map induces the identity on homology groups. This conjugation is by construction trivial.
\end{Proof}

\noindent Note that the above two arguments work in parallel for the two cases of groups we are considering here, we only need Corollary~\ref{c:wittun} and Proposition~\ref{p:wittgln}.

\section{Acyclicity of the orbit chain complex}

We have already shown that the bottom row of $E^1$ consists of a chain complex which is isomorphic to $(\Z B_p)_p$, the complex of free abelian groups over the $G$-orbits of $X^{p+1}$. The main objective of this section will be to show that the bottom row of the lower triangle in $E^2$ is always zero, we will, however, need a stronger result for the induction step.

For the further course of the proof, we note that the dimension filtration $\cF^k$ induces a filtration on $E^1$. Let 
\[
	\cF^k_j = \cF^k \bigl(E^1_{*,j}\bigr)
\]
denote the associated filtered `row'-complex, which is a chain complex with respect to the restriction of the differential $d^1$, which, as we know, is induced by $\id_{F(G)} \otimes_G d_C$.

The required stronger result is the following
\begin{Lemma}
	$\cF^k_0$ is $k$-acyclic with augmentation $H_0(G,C_0)$ for any $1\leq k\leq n$. For $k=n$ the filtration is trivial, as every simplex spans at most an $(n+1)$-dimensional subspace. This especially yields
\[
	E^2_{p,0}=0\quad\text{for}\quad 1\leq p\leq n.
\] \label{l:zeroesbottom}
\end{Lemma}
\begin{Proof}
We use equation \eqref{e:bottomrow}. The dimension-filtered complex can then easily be written as
\begin{equation}
	(\cF^k_0)_p \cong \cF^k (\Z B_p) = \Z (\cF^k B_p). \label{e:dimfiltorbcplx}
\end{equation}
 Every element in $\cF^k_0$ can hence be expressed as a linear combination of elements of the form $[c]$, where $c$ is any simplex satisfying $|[c]|=|c|\leq k+1$.

Now consider the short exact sequence of chain complexes
\[
	0 \rightarrow \cF^{k-1}_0 \stackrel{\iota}{\rightarrow} \cF^k_0 \rightarrow \cF^k_0/\cF^{k-1}_0 \rightarrow 0,
\]
where $\iota$ is the canonical inclusion. We get a long exact homology sequence
\[
	\cdots \rightarrow H_l(\cF^{k-1}_0) \stackrel{\iota_*}{\rightarrow} H_l(\cF^k_0) \rightarrow H_l(\cF^k_0/\cF^{k-1}_0) \rightarrow H_{l-1}(\cF^{k-1}_0) \stackrel{\iota_*}{\rightarrow} H_{l-1}(\cF^k_0) \rightarrow \cdots
\]

\noindent Our first goal will be to show that $\iota_*$ is the zero map.

\subsection{Breaking up the long exact homology sequence}

Of course, it is enough to show that the inclusion $\iota : \cF^{k-1}_0\rightarrow \cF^k_0$ maps all generating cycles $[c]$ onto boundaries. We note that for any element $[c] \in \cF^{k-1}_0$, its dimension $|[c]|$ is at most $k$, by definition of the filtration, so
\[
\im\iota = \langle [c]\in\cF^k_0 : |[c]| \leq k \rangle.
\]
This means that, in the general linear case, for any $[c]\in\im\iota$, we can pick a vector $v\in X$ with $v \not\in \spn c$ , since $|[c]|\leq k \leq n$. In the unitary case, we can even pick $v\in (\spn c)^\perp$. Take any $[c]\in\im\iota\backslash (\cF^k_0)^{\phantom k}_0$, then we have $c\not\in C_0$. We calculate the differential of $[v \msharp c]\in\cF^k_0$, the join of $c$ with $v$ (see Definition \ref{d:msharp}):
\[
	d^1[v \msharp c] = [d_C(v \msharp c)] \stackrel{c\not\in C_0}{=} [c - (v \msharp (d_C c)] = [c] - [v \msharp d_C c].
\]
We know that $[d_C c] = 0$, since $[c]$ is a cycle. Therefore the simplices in $d_Cc$ cancel out under the action of $G$. Since $\spn(c)$ is contained in an $n$-dimensional hyperplane, this cancellation can be realized by invertible or, respectively, unitary maps fixing $v$ by Corollary~\ref{c:wittun} and Proposition~\ref{p:wittgln}. So $[v\msharp d_Cc] = 0$, hence $[c]$ is a boundary. So $\iota_*:H_l(\cF^{k-1}_0)\rightarrow H_l(\cF^k_0)$ is the zero homomorphism for $l>0$.

It follows that, for $l>1$, the long homology sequence breaks up into short exact pieces
\[
	0 \rightarrow H_l(\cF^k_0) \rightarrow H_l(\cF^k_0/\cF^{k-1}_0) \rightarrow H_{l-1}(\cF^{k-1}_0) \rightarrow 0.
\]
For $l=1$ we have an injective map $H_1(\cF^k_0) \rightarrow H_1(\cF^k_0/\cF^{k-1}_0)$, which is all we need.
For a proof of the $k$-acyclicity of $\cF^k_0$, it is now enough to prove the $k$-acyclicity of $\cF^k_0/\cF^{k-1}_0$ with respect to the differential $\partial$ induced by $d^1$.

\subsection{Acyclicity of the factor complex}\label{ss:acyclfact}

We need to show that $\cF^k_0/\cF^{k-1}_0$ is $k$-acyclic. Since $(\cF^k_0/\cF^{k-1}_0)_l=0$ for $l<k$, these homology groups are zero as well. The only thing we need to prove is that $H_k(\cF^k_0/\cF^{k-1}_0)=0$. 

It is clear that each element in $(\cF^k_0/\cF^{k-1}_0)_k$ is an equivalence class of a sum of simplices consisting of $(k+1)$ linearly independent vectors. The idea is to pick such a simplex which is a cycle and show that it is already a boundary.

Let again $[c]$ denote the equivalence class of such a cycle $c=(v_0,\ldots,v_k)$. Although this is a slight abuse of notation, it will only be used in this part and will never conflict with the definition above. 

\ruleboxstart
\minisec{The general linear case}
Consider the following calculation:
\begin{align*}
	\partial[(v_0+v_1,v_0,v_1,v_2,\ldots,v_k)] =& [d_C(v_0+v_1,v_0,v_1,v_2,\ldots,v_k)] \\
	=& [(v_0,v_1,\ldots,v_k)] \\
	&- [(v_0+v_1,v_1,\ldots,v_k)] + [(v_0+v_1,v_0,v_2,\ldots,v_k)]\\
	&- [(v_0+v_1,v_0,v_1,v_3,\ldots,v_k)] + \ldots \\
	&\pm[(v_0+v_1,v_0,v_1,\ldots,v_{k-1})]
\end{align*}
The terms in the last two lines all vanish, since the vectors are not linearly independent. We get
\[
\partial[(v_0+v_1,v_0,v_1,v_2,\ldots,v_k)]= [c] - [(v_0+v_1,v_1,\ldots,v_k)] + [(v_0+v_1,v_0,v_2,\ldots,v_k)].
\]
Since any $(k+1)$-tuple of linearly independent vectors can be mapped onto any other such tuple by a linear map in $\Gl_{n+1}\F$, the last two equivalence classes are the same and cancel out. Hence $c$ is a boundary, the factor complex is acyclic.

\minisec{The unitary case} This is a point in the proof where the unitary case is much more difficult. Let
\[
	B=\partial((\cF^k_0/\cF^{k-1}_0)_{k+1})\subseteq (\cF^k_0/\cF^{k-1}_0)_k
\]
denote the submodule of boundaries. We also need an auxiliary definition:
\begin{Definition}
	For an $l$-simplex $c=(v_0,\ldots,v_l)$ we consider the space
\[
	\Big(\sum_{j=1}^l(v_j-v_0)\K\Big)^\perp \cap \spn(v_0,\ldots,v_l).
\]
Since this space is obviously one-dimensional, we denote it by $v\K$. Such a spanning $v$ fulfils the relation
\[
	\langle v, v_0 \rangle = \langle v, v_j \rangle = r \qquad\forall\; 1\leq j\leq l.
\]
We now pick a $v$ such that $\langle v,v\rangle=1$ and $r\in\R$ with $r>0$\footnote{$r$ being positive is not important for the course of this proof. This is, however, necessary for $\cC(c)$ to be well-defined.}. This clearly determines $v\in S^n$ uniquely, this unique $v$ is called the \emph{circumcentre $\cC(c):=v$} of the simplex $c$.
\end{Definition}

Our goal is to show that any cycle $[c]$ is a boundary, i.e.~$[c]\in B$. To do so, we modify $[c]$ step by step by boundaries until we reach a boundary.

To be more explicit, let $z_0=\cC(c)$, then we know that:
\[
	\partial [z_0 \msharp c] = [c] - [z_0 \msharp d_Cc].
\]
So we have
\[
	[c] \equiv_B [z_0 \msharp d_Cc],
\]
where the right side is a sum of equivalence classes of $(k-1)$-simplices joint with $z_0$. If we pick one summand, call it $[z_0\msharp(w_1,\ldots,w_k)]$, we know that
\[
	\langle z_0,w_1\rangle = \langle z_0,w_j \rangle\qquad \forall\; 2\leq j \leq k
\]
by construction of the circumcentre $z_0$. If we continue with choosing
\[
	z_1=\cC(w_1,\ldots,w_k), 
\]
we see that
\begin{equation*}
	\begin{split}
		\partial[z_0\msharp(z_1 \msharp(w_1,\ldots,w_k))] = &\underbrace{[z_1 \msharp (w_1,\ldots,w_k)]}_{=0,\text{ since } \dim<k+1} - [z_0\msharp(w_1,\ldots,w_k)] \\&+ [z_0\msharp(z_1\msharp d_C(w_1,\ldots,w_k))],
	\end{split}
\end{equation*}
and we arrive at
\[
	[z_0\msharp(w_1,\ldots,w_k)] \equiv_B [z_0 \msharp ( z_1 \msharp d_C(w_1,\ldots w_k))].
\]
Again, we note that $d_C(w_1,\ldots,w_k)$ is a sum of $(k-2)$-simplices, of which every summand $(x_2,\ldots,x_k)$ has the property
\[
	\langle z_0,x_2\rangle = \langle z_0,x_j \rangle,\;\langle z_1,x_2 \rangle = \langle z_1,x_j \rangle \qquad \forall \;3\leq j \leq k
\]
If we continue this procedure, always replacing entries with the circumcentres of the `tail', we finally arrive at
\begin{equation*}
	\begin{split}
	[c] &\equiv_B \sum_{j\in J} \big[z_0^j \msharp( z_1^j \msharp (\cdots(z_{k-2}^j \msharp(u_{k-1}^j,u_k^j))\cdots))\big] \\
		&= \sum_{j\in J} \big[(z_0^j,\ldots,z_{k-2}^j,u_{k-1}^j,u_k^j)\big] \\
	\end{split}
\end{equation*}
for a finite index set $J$ with
\[
	\langle z_i^j, u_{k-1}^j \rangle = \langle z_i^j,u_k^j \rangle\qquad \forall \;0\leq i\leq k-2,\quad \forall j\in J.
\]
If, again, we pick $z_{k-1}^j=\cC((u_{k-1}^j,u_k^j))$ and modify by a boundary as before, we have
\[
	[c] \equiv_B \sum_{j\in J} (-1)^k\big( [z_0^j,\ldots,z_{k-1}^j,u_k^j] - [z_0^j,\ldots,z_{k-1}^j,u_{k-1}^j]\big).
\]
Because of the above inner product relations and the construction of the circumcentre $z_{k-1}^j$, Corollary~\ref{c:wittun} implies that the last two simplices are congruent modulo $G$. This implies $[c] \equiv_B 0$, hence $[c]\in B$.
\ruleboxend

\paragraph{Summary} Putting together all the information, we get that $\cF^k_0$ is $k$-acyclic for any $1\leq k \leq n$ in both cases.
\end{Proof}

\section{Induction}

The rest of the proof will be by induction on $n$. It will also work completely in parallel for both the unitary and the general linear case.

\subsection{Base case}

The base case is not very difficult to see --- just remark that for $n=1$, the conditions of the triangle lemma (Lemma~\ref{l:triangle}) are already fulfilled. The `triangle' in this case is rather small, all that is required is that $E^2_{1,0}$ vanishes, which has been proven in Lemma~\ref{l:zeroesbottom}.
\begin{figure}[htb]
	\input{basecase} \caption{The base case situation} \label{f:basecase}
\end{figure}
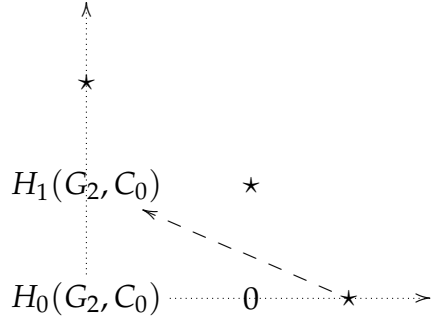

\noindent So by the triangle lemma, the edge map
\begin{align*}
	\pi_0: H_0(G_2,C_0) = E^2_{0,0} &\rightarrow E^\infty_{0,0} = H_0(G_2,C) \\
\intertext{is an isomorphism induced by}
	f \otimes_{G_2} zc &\mapsto f \otimes_{G_2} zc.
\end{align*}
By Proposition~\ref{p:homologyofg} there is an isomorphism
\begin{align*}
	H_q(G_2,C) &\rightarrow H_q(G_2) \\
\intertext{induced by}
	f \otimes_{G_2} zc &\mapsto f \otimes_{G_2} \varepsilon_{F(G_2)}(zc) = f \otimes_{G_2} z
\end{align*}
for all $q\in\Z$. In Lemma~\ref{l:specseqstruct} we had constructed isomorphisms
	\begin{align*}
		H_q(G_1) &\rightarrow H_q(G_2,C_0)=E^1_{0,q}=E^2_{0,q}\\
		\intertext{induced by}
		f \otimes_{G_1} z &\mapsto (i^{G_1}_{G_2})f \otimes_{G_2} z c.
	\end{align*}
Concatenating all these isomorphisms, we get an isomorphism
\begin{align*}
	H_0(G_1) &\rightarrow H_0(G_2)\\
\intertext{induced by}
	f \otimes_{G_1} z &\mapsto (i^{G_1}_{G_2})f \otimes_{G_2} z,
\end{align*}
so this is exactly the isomorphism $(i^{G_1}_{G_2})_\stern$, which proves the first part of the theorem.
 The case for $H_1$ is not much more difficult, here the edge map
\[
	\pi_1: H_1(G_2,C_0) = E^2_{0,1} \rightarrow E^\infty_{0,1} = H_1(G_2,C)
\]
is surjective, hence the map $(i^{G_1}_{G_2})_\stern$ is just surjective.

\subsection{Induction step}

Now, assume that the theorem is proven for all numbers up to $n$. As before, we consider the inclusion
\[
	(i^{G_n}_G) : G_n \rightarrow G=G_{n+1}
\]
and the spectral sequence $E$ associated to $G$ and $C(X)$, which we have constructed before. We want to apply the triangle lemma. We already know that, in $E^2$, the leftmost column is isomorphic to $H_\stern(G_n)$ and the bottom row consists of zeroes.

The rest of the proof will be concerned with filling the triangle with zeroes, that is
\begin{Lemma} We have
	\[
	E^2_{p,q}=0,\qquad\forall p+q\leq n,\quad p>0.
	\]\label{l:lowertriangle}
\end{Lemma}
\begin{Proof} We need to prove that the subcomplexes $E_{\star,q}^1$ are $(n-q)$-acyclic for all $1\leq q \leq n-1$. The dimension-filtered complex $\cF^{n-q}_q$ differs from the unfiltered complex only in the boundaries in $(\cF^{n-q}_q)_q$. There are fewer of them, namely only those of simplices with dimension less than or equal to $q+1$. If this complex is $(n-q)$-acyclic, surely the whole complex $E^1_{\stern,q}$ is $(n-q)$-acyclic. We show that $\cF^{n-q}_q$ is $(n-q)$-acyclic for all $1\leq q\leq n-1$.

So, let $1\leq q \leq n-1$ be arbitrary. Then by Lemma~\ref{l:specseqstruct}, we know that
\[
	(\cF^{n-q}_q)_p \cong \bigoplus_{c\in \cF^{n-q}B_p} H_q(G_{(n+1)-|c|}).
\]

\noindent By definition of the filtration, for $c\in \cF^{n-q}B_p$ we have $1\leq |c| \leq n+1-q$ and so
\[
	H_q(G_{(n+1)-|c|}) \cong 
	\left\{ \begin{array}{ll} 
		H_q(G_q) & \text{ for } |c|=n+1-q \\
		H_q(G_{q+1}) & \text{ for } |c|<n+1-q.
\end{array}\right.
\]
The lower isomorphism is a part of the induction hypothesis, since 
\[
	H_q(G_{q+1})\cong H_q(G_r) \quad\text{for}\quad q+1\leq r\leq n.
\]
Hence we know $\cF^{n-q}_q$ is isomorphic to a direct sum of the groups $H_q(G_q)$ and \\$H_q(G_{q+1})$, depending on the dimension of the corresponding simplex. We have
\[
	\cF^{n-q}_q \cong \bigoplus_{c\in \cF^{n-q}B_p} \left\{\begin{array}{cc} H_q(G_q) & |c|=n+1-q \\ H_q(G_{q+1}) & |c|<n+1-q \end{array}\right..
\]
The other part of the induction hypothesis gives us an epimorphism
\[
	i_\stern: H_q(G_q)\twoheadrightarrow H_q(G_{q+1}),
\]
hence an exact sequence
\[
	0 \rightarrow K_q \rightarrow H_q(G_q)\stackrel{i_\stern}{\rightarrow} H_q(G_{q+1}) \rightarrow 0.
\]
Let $K_q\leq H_q(G_q)$ denote the kernel of this map. We already know that
\[
	(\cF^{n-q}_0)_p \cong \Z(\cF^{n-q} B_p),
\]
the free abelian groups over $G$-orbits, by equation \eqref{e:dimfiltorbcplx}. We can write down the following map
\begin{align*}
	H_q(G_q)\otimes \cF^{n-q}_0 &\rightarrow \cF^{n-q}_q \\
	h \otimes \sum_c z_c [c] & \mapsto \sum_c \left\{\begin{array}{cc} h & |c|=n+1-q \\ i_\stern(h) & |c|<n+1-q \end{array}\right.
\end{align*}
This map is by construction surjective, the kernel is $K_q\otimes \cF^{n-(q+1)}_0$. We get a short exact sequence of chain complexes
\[
	0 \rightarrow K_q\otimes \cF^{n-(q+1)}_0 \stackrel{I}{\rightarrow} H_q(G_q)\otimes \cF^{n-q}_0 \rightarrow \cF^{n-q}_q  \rightarrow 0,
\]
where $I$ is the inclusion map.
Associated to this short exact sequence of chain complexes, we get a long exact homology sequence:
\begin{equation*}
	\begin{split}
	\cdots &\rightarrow H_l(K_q\otimes \cF^{n-(q+1)}_0) \rightarrow H_l(H_q(G_q)\otimes \cF^{n-q}_0) \rightarrow H_l(\cF^{n-q}_q) \rightarrow \\ &\rightarrow H_{l-1}(K_q\otimes \cF^{n-(q+1)}_0) \rightarrow \cdots
	\end{split}
\end{equation*}
We want to prove that $H_l(\cF^{n-q}_q)=0$ for $1\leq l \leq n-q$. We first examine the other terms.

\noindent The chain complex $\cF^{n-(q+1)}_0$  consists of free abelian groups, we can hence apply the universal coefficient theorem:
\[
	0 \rightarrow K_q\otimes H_l(\cF^{n-(q+1)}_0) \rightarrow H_l(K_q\otimes \cF^{n-(q+1)}_0) \rightarrow \Tor(H_{l-1}(\cF^{n-(q+1)}_0),K_q) \rightarrow 0
\]
By Lemma~\ref{l:zeroesbottom}, we have
\[
	H_l(\cF^{n-(q+1)}_0)=0\quad\text{for}\quad 1\leq l \leq (n-(q+1)),
\]
so the left terms in these short sequences vanish for $1\leq l \leq n-(q+1)$. Then for $l\neq 1$ the $\Tor$-terms vanish as well. In addition, $H_0(\cF^{n-(q+1)})$ is free, since it is a submodule of $\cF^{n-(q+1)}_0$, which is free, the $\Tor$-term also vanishes for $l=1$. Exactly the same argument works for the $H_l(H_q(G_q)\otimes \cF^{n-q}_0)$-terms.

So, taking a look at the long exact sequence above, we see that 
\[
	H_l(\cF^{n-q}_q)=0\quad\text{for}\quad 2\leq l \leq n-q.
\]
The only remaining problem is $H_1$. We look at the corresponding part of the long homology sequence:
\[
	\cdots\rightarrow 0 \rightarrow H_1(\cF^{n-q}_q) \rightarrow H_0(K_q\otimes \cF^{n-(q+1)}_0) \rightarrow H_0(H_q(G_q)\otimes \cF^{n-q}_0) \rightarrow \cdots
\]
We again apply the universal coefficient theorem. As all $\Tor$-terms vanish, we arrive at:
\[
	\cdots\rightarrow 0 \rightarrow H_1(\cF^{n-q}_q) \rightarrow K_q\otimes H_0(\cF^{n-(q+1)}_0) \stackrel{I_\stern}{\rightarrow} H_q(G_q)\otimes H_0(\cF^{n-q)}_0) \rightarrow \cdots
\]
The problem is to show that $I_\stern$ is injective. We write down the short exact sequence
\[
	0 \rightarrow \cF^{n-(q+1)}_0 \rightarrow \cF^{n-q}_0 \rightarrow \cF^{n-q}_0 / \cF^{n-(q+1)}_0 \rightarrow 0
\]
In section \ref{ss:acyclfact} we already proved that $H_1(\cF^{n-q}_0 / \cF^{n-(q+1)}_0)$ vanishes, since we have $n-q\geq 1$. So, by the long exact sequence, the inclusion-induced map
\[
	H_0(\cF^{n-(q+1)}_0) \rightarrow H_0(\cF^{n-q}_0)
\]
is injective. Then $I_\stern$ is also injective as it is the tensor product of two injective maps. Hence $\cF^{n-q}_q$ is $(n-q)$-acyclic for all $1\leq q\leq n-1$, which we wanted to prove.
\end{Proof}
\begin{figure}[bht]
\input{triangle2} \caption{Induction step: situation} \label{f:triangle2}
\end{figure}
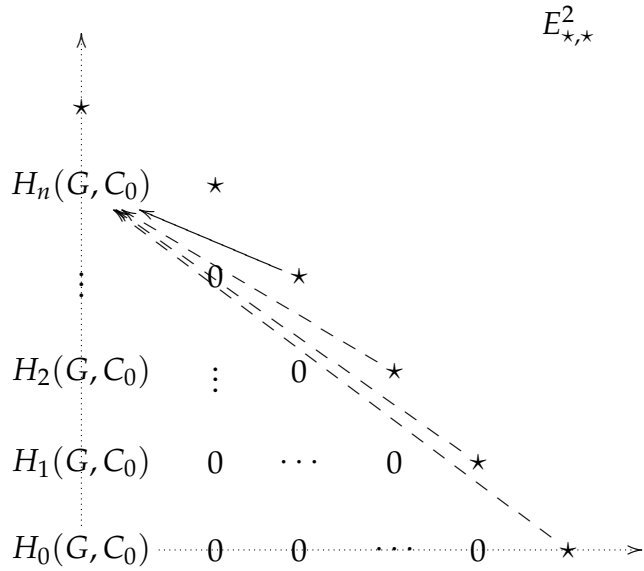
\paragraph{Construction of the isomorphisms} We just proved Lemma~\ref{l:lowertriangle}, that means the lower triangle of $E^2$ is filled with zeroes, save the leftmost column, which consists of the groups $H_q(G,C_0)$. To make that explicit, we have the following situation
	\[
		E^2_{p,q} = \left\{\begin{array}{ccc}
			H_q(G,C_0) & p=0,& 0\leq q \leq n \\
			0 & p>0, &p+q\leq n
		\end{array}\right.,
	\]
which is exactly the prerequisite for the triangle lemma, see figure \ref{f:triangle2}.

We apply the triangle lemma (Lemma~\ref{l:triangle}). Hence the edge homomorphism
	\begin{align*}
		\pi_q: H_q(G,C_0) = E^2_{0,q} &\twoheadrightarrow E^\infty_{0,q} = H_q(G,C) \\
		\intertext{induced by}
		f \otimes_G z c &\mapsto f\otimes_G z c
	\end{align*}
is surjective for $1\leq q \leq n$. It is an isomorphism for $1\leq q < n$.

\noindent As in the base case, we have isomorphisms
\begin{align*}
	H_q(G,C) &\rightarrow H_q(G) \\
\intertext{induced by}
	f \otimes_G z c &\mapsto f \otimes_G \varepsilon_{F(G)}(zc) = f\otimes_G z
\end{align*}
by Proposition~\ref{p:homologyofg} and isomorphisms
\begin{align*}
		H_q(G_n) &\rightarrow H_q(G,C_0)=E^1_{0,q}=E^2_{0,q}\\
		\intertext{induced by}
		f \otimes_{G_n} z &\mapsto (i^{G_n}_{G})f \otimes_{G} z c
\end{align*}
by Lemma~\ref{l:specseqstruct}.

\noindent Putting it all together, the concatenation of these three maps yields, as above, that
\[
	(i^{G_n}_G)_\stern: H_q(G_n)\rightarrow H_q(G)
\]
is an isomorphism for $q<n$ and an epimorphism for $q=n$. This finishes the proof of the theorem.\hfill$\Box$

\paragraph{Remark} In Sah's original proof for the general linear case, a sharper stability range is achieved by proving a stronger version of Lemma~\ref{l:zeroesbottom}, namely that $\cF^k_0$ is $(k+1)$-acyclic for $1\leq k \leq n$, if $n\geq 2$, see \cite[Theorem~B.7]{Sah:HcL:86}. Then the triangle lemma can be applied to a larger triangle, therefore yielding a stronger stability result. But including this stronger version would have disturbed the symmetry in our arguments, therefore we excluded this part.

\section{Consequences}

To finish this thesis, we want to list some easy consequences of the Theorems \ref{t:unstability} and \ref{t:glnstability}.

Throughout this section, let $M$ be an abelian group, considered to be a trivial $G$-module for any group $G$. Then, of course, $H_k(-,M)$ is a covariant functor. The argument is completely analogous to the case $M=\Z$. Hence we can talk about homology stability theorems for any trivial coefficient group.

\begin{BigTheorem}[Stability of homology with coefficients for unitary groups]
For $n\geq 0$, let $U_n\F$ be the standard unitary group over $\F$, where $\F \in \{\R,\C,\bH\}$. The canonical inclusion
	\[
	i^{U_n\F}_{U_{n+1}\F}: U_n \F \rightarrow U_{n+1}\F
	\]
	induces a homology isomorphism
	\[
	(i^{U_n\F}_{U_{n+1}\F})_*: H_q(U_n \F,M) \rightarrow H_q(U_{n+1}\F,M)
	\]
	for $q< n$. The induced map is surjective for $q\leq n$.
\end{BigTheorem}

\begin{BigTheorem}[Stability of homology with coefficients for general linear groups]
For $n\geq 0$, let $\Gl_n\F$ be the $n$-dimensional general linear group over $\F$, where $\F$ is any skew-field with infinite centre. The canonical inclusion
	\[
	i^{\Gl_n\F}_{\Gl_{n+1}\F}: \Gl_n \F \rightarrow \Gl_{n+1}\F
	\]
	induces a homology isomorphism
	\[
	(i^{\Gl_n\F}_{\Gl_{n+1}\F})_*: H_q(\Gl_n \F,M) \rightarrow H_q(\Gl_{n+1}\F,M)
	\]
	for $q<n$. The induced map is surjective for $q\leq n$.
\end{BigTheorem}

\begin{Proof}
	We use the universal coefficient theorem and abbreviate $G_n=U_n\F$ or $G_n=\Gl_n\F$, respectively. For $q\leq n$ we have
	\[
	\xymatrix{
		0\ar[r] & M\otimes H_q(G_n)\ar[r]\ar[d]^{i_\stern} & H_q(G_n,M)\ar[r]\ar[d]^{i_\stern} & \Tor(H_{q-1}(G_n),M)\ar[r]\ar[d]^{i_\stern} & 0 \\
		0\ar[r] & M\otimes H_q(G_{n+1})\ar[r] & H_q(G_{n+1},M)\ar[r] & \Tor(H_{q-1}(G_{n+1}),M)\ar[r] & 0
	}
	\]
	By Theorems~\ref{t:unstability} and \ref{t:glnstability}, the maps $i_\stern:H_{q-1}(G_n)\rightarrow H_{q-1}(G_{n+1})$ are all isomorphisms under the hypothesis. Then the induced maps on $\Tor$ are isomorphisms as well. The maps $i_\stern:H_q(G_n)\rightarrow H_q(G_{n+1})$ are isomorphisms or epimorphisms, depending on $q$. Apply the five-lemma and conclude the two theorems.
\end{Proof}

\noindent Similarly, according theorems for group cohomology can now be proven.

\begin{BigTheorem}[Cohomological stability for unitary groups]
	For $n\geq 0$, let $U_n\F$ be the standard unitary group over $\F$, where $\F \in \{\R,\C,\bH\}$. The canonical inclusion
	\[
	i^{U_n\F}_{U_{n+1}\F}: U_n \F \rightarrow U_{n+1}\F
	\]
	induces a homology isomorphism
	\[
		(i^{U_n\F}_{U_{n+1}\F})^\stern: H^q(U_{n+1} \F,M) \rightarrow H^q(U_n\F,M)
	\]
	for $q<n$. The induced map is injective for $q\leq n$.
\end{BigTheorem}

\begin{BigTheorem}[Cohomological stability for general linear groups]
	For $n\geq 0$, let $\Gl_n\F$ be the $n$-dimensional general linear group over $\F$, where $\F$ is any skew-field with infinite centre. The canonical inclusion
	\[
	i^{\Gl_n\F}_{\Gl_{n+1}\F}: \Gl_n \F \rightarrow \Gl_{n+1}\F
	\]
	induces a homology isomorphism
	\[
		(i^{\Gl_n\F}_{\Gl_{n+1}\F})^\stern: H^q(\Gl_{n+1} \F) \rightarrow H^q(\Gl_n\F)
	\]
	for $q<n$. The induced map is injective for $q\leq n$.
\end{BigTheorem}
\begin{Proof}
	We again write down the short exact sequences from the universal coefficient theorem for any $q\leq n$, abbreviating $G_n=U_n\F$ or $G_n=\Gl_n\F$.
	\[
	\xymatrix{
		0 & \Ext(H_{q-1}(G_{n+1}),M)\ar[l]\ar[d]^{i^\stern} & H^q(G_{n+1},M)\ar[l]\ar[d]^{i^\stern} & \Hom(H_q(G_{n+1}),M)\ar[l]\ar[d]^{i^\stern} & 0\ar[l] \\
		0 & \Ext(H_{q-1}(G_n),M)\ar[l] & H^q(G_n,M)\ar[l] & \Hom(H_q(G_n),M)\ar[l] & 0\ar[l]
	}
	\]
	By the Theorems~\ref{t:unstability} and \ref{t:glnstability}, the map
	\[
		i_\stern: H_{q-1}(G_n)\rightarrow H_{q-1}(G_{n+1})
	\]
	is always an isomorphism under the hypothesis, hence so is $i^\stern$ on $\Ext$. By the same token
	\[
		i_\stern: H_q(G_n)\rightarrow H_q(G_{n+1})
	\]
	is an isomorphism for $q<n$ and an epimorphism for $q=n$. The induced map $i^\stern$ on $\Hom$ is hence an isomorphism for $q<n$ and a monomorphism for $q=n$.
	
	We get the desired result by applying the five-lemma.
\end{Proof}
Since the group homology functor commutes with direct limits (see for example \cite[page 195]{Bro:CoG:82}), we also have the following easy consequences:
\begin{BigTheorem}[Homology of \textnormal{$\Gl\F$}] For any skew-field with infinite centre $\F$, the canonical inclusion $\Gl_n\F \rightarrow \Gl\F$ induces a homology isomorphism
\[
	H_q(\Gl_n\F) \cong H_q(\Gl\F)
\]
for $q<n$.
\end{BigTheorem}

\begin{BigTheorem}[Homology of \textnormal{$U\F$}] For $\F \in \{ \R,\C,\bH \}$, the canonical inclusion $U_n\F \rightarrow U\F$ induces a homology isomorphism
\[
	H_q(U_n\F) \cong H_q(U\F)
\]
for $q<n$.
\end{BigTheorem}

\nocite{Ros:AKa:94}
\nocite{Hut:naM:90}
\bibliographystyle{alpha}
\bibliography{biblio}
\end{document}

%% file: firstpages.tex
\[
\begin{array}{c|c}
\begin{xy}
	\xymatrix@=1.2em{
	&&&&&E^0_{\star,\star}\\
	\ar[d]\\
	E^0_{0,3}\ar[d]&\ar[d]\\
	E^0_{0,2}\ar[d]&E^0_{1,2}\ar[d]&\ar[d]\\
	E^0_{0,1}\ar[d]&E^0_{1,1}\ar[d]&E^0_{2,1}\ar[d]&\ar[d]\\
	E^0_{0,0}\ar@.[uuuuu]\ar@.[rrrrr]&E^0_{1,0}&E^0_{2,0}&E^0_{3,0}&&
}
\end{xy}&
\begin{xy}
\xymatrix@=1.2em{
	&&&&&E^1_{\star,\star}\\
	\\
	E^1_{0,3}&\ar[l]\\
	E^1_{0,2}&E^1_{1,2}\ar[l]&\ar[l]\\
	E^1_{0,1}&E^1_{1,1}\ar[l]&E^1_{2,1}\ar[l]&\ar[l]\\
	E^1_{0,0}\ar@.[uuuuu]\ar@.[rrrrr]&E^1_{1,0}\ar[l]&E^1_{2,0}\ar[l]&E^1_{3,0}\ar[l]&\ar[l]&
}
\end{xy}\\\hline
\begin{xy}
\xymatrix@=1.2em{
	&&&&&E^2_{\star,\star}\\
	\\
	E^2_{0,3}\\
	E^2_{0,2}&E^2_{1,2}&\ar[llu]\\
	E^2_{0,1}&E^2_{1,1}&E^2_{2,1}\ar[llu]&\ar[llu]\\
	E^2_{0,0}\ar@.[uuuuu]\ar@.[rrrrr]&E^2_{1,0}&E^2_{2,0}\ar[llu]&E^2_{3,0}\ar[llu]&\ar[llu]&
}
\end{xy}&

\end{array}
\]

%% file: convergence.tex
\[	
	\xymatrix@=1.2em{
	0&&\\
	&&E^3_{0,2}\\
	&&E^3_{0,1}& \ar[llluu]E^3_{1,1}\\
	&&E^3_{0,0}\ar@.[uuu]\ar@.[rrr] & E^3_{1,0} & E^3_{2,0}&\\
	&&&&&&0\ar[llluu]
	}
\]

%% file: edgehom.tex
\[
	\xymatrix@=1.2em{
	&&\\
	&&\ar[l]\ar[llu]E_{0,2}\\
	&&E_{0,1}&E_{1,1}\\
	&&E_{0,0}\ar@.[uuu]\ar@.[rrr] & E_{1,0} & E_{2,0}&\\
	&&&&&\ar[llu]\\
	&&&&&&\ar[llluu]
	}
\]

%% file: triangle.tex
\[
\begin{xy}
\xymatrix@=1em{
	&&&&&E^2_{\star,\star}\\
	\star\\
	G_n&\star\\
	\vdots&0&\star\ar[llu]\\
	G_2&\vdots&0&\star\ar@{-->}[llluu]&\\
	G_1&0&\cdots&0&\star\ar@{-->}[lllluuu]\\
	G_0\ar@.[uuuuuu]\ar@.[rrrrrr]&0&0&\cdots&0&\star\ar@{-->}[llllluuuu]&
}
\end{xy}
\]

%% file: basecase.tex
\[
\begin{xy}
\xymatrix{
	\\
	\star\\
	H_1(G_2,C_0)&\star\\
	H_0(G_2,C_0)\ar@.[uuu]\ar@.[rrr]&0&\ar@{-->}[llu]\star&
}
\end{xy}
\]

%% file: triangle2.tex
\[
\begin{xy}
\xymatrix@=1.2em{
	&&&&&E^2_{\star,\star}\\
	\star\\
	H_n(G,C_0)&\star\\
	\vdots&0&\star\ar[llu]\\
	H_2(G,C_0)&\vdots&0&\star\ar@{-->}[llluu]&\\
	H_1(G,C_0)&0&\cdots&0&\star\ar@{-->}[lllluuu]\\
	H_0(G,C_0)\ar@.[uuuuuu]\ar@.[rrrrrr]&0&0&\cdots&0&\star\ar@{-->}[llllluuuu]&
}
\end{xy}
\]